%% file: agt-2-33.tex
\documentclass{gtart}
\input agtout

\lognumber{33}
\volumenumber{2}
\volumeyear{2002}
\papernumber{33}
\published{3 October 2002}
\pagenumbers{791}{824}
\received{6 March 2002}
\accepted{4 September 2002}

\usepackage{amsmath, amssymb, graphicx}
\theoremstyle{plain}
\newtheorem{thm}{Theorem}[section]
\newtheorem{lem}[thm]{Lemma}
\newtheorem{cor}[thm]{Corollary}

\newcommand{\Diff}{\text{Diff}}

\newcommand{\Ortho}{\text{O}}
\newcommand{\Aut}{\text{Aut}}
\newcommand{\Symp}{\text{Sp}}
\newcommand{\inv}[1]{\overline{#1}\;}
\newcommand{\Prod}[3]{\prod_{#1}^{#2} #3}
\newcommand{\Ggeq}{\underset{G_g}{\sim}}
\newcommand{\Ccirc}{\circ\circ\circ}
\newcommand{\Cbullet}{\bullet\bullet\bullet}
\newcommand{\repre}[1]{\overline{\overline{#1}}}
\newfont{\ttnonit}{cmtt10}
\theoremstyle{definition}

\begin{document}
\title{On diffeomorphisms over surfaces trivially\\embedded in 
the 4-sphere}
\author{Susumu Hirose}
\address{Department of Mathematics, 
Faculty of Science and Engineering\\Saga University, 
Saga, 840-8502 Japan.}
\email{hirose@ms.saga-u.ac.jp}
\begin{abstract} 
A surface in the 4-sphere is {\em trivially\/} embedded, 
if it bounds a 3-dimensional handle body in the 4-sphere. 
For a surface trivially embedded in the 4-sphere, a diffeomorphism 
over this surface is extensible if and only if this preserves 
the Rokhlin quadratic form of this embedded surface. 
\end{abstract}
\asciiabstract{ 
A surface in the 4-sphere is trivially embedded, 
if it bounds a 3-dimensional handle body in the 4-sphere. 
For a surface trivially embedded in the 4-sphere, a diffeomorphism 
over this surface is extensible if and only if this preserves 
the Rokhlin quadratic form of this embedded surface.}
\maketitle
\cl{\small\it This paper is 
dedicated to Professor Mitsuyoshi Kato on his 60th birthday.}

\section{Introduction}
We denote the closed oriented surface of genus $g$ by $\Sigma_g$, 
the mapping class group of $\Sigma_g$ by ${\cal M}_g$. 
Let $\phi\co \Sigma_g \to S^4$ be an embedding, and 
$K$ be its image. 
We call $(S^4, K)$ a {\it $\Sigma_g$-knot \/}. 
Two $\Sigma_g$-knots $(S^4, K)$ and $(S^4, K')$ are 
{\it equivalent\/} if there is a diffeomorphism of $S^4$ 
which brings $K$ to $K'$. 
{\it A 3-dimensional handlebody\/} $H_g$ is an oriented 3-manifold 
which is constructed from a 3-ball with attaching $g$ 1-handles. 
Any embeddings of $H_g$ into $S^4$ are isotopic each other. 
Therefore, $(S^4, \partial H_g)$ is unique up to equivalence. 
We call this $\Sigma_g$-knot $(S^4, \partial H_g)$ 
{\it a trivial $\Sigma_g$-knot\/} and 
denote this by $(S^4, \Sigma_g)$. 
For a $\Sigma_g$-knot $(S^4, K)$, we define the following group, 
\begin{align*}
{\cal E}(S^4, K) 
= \left\{ \phi \in \pi_0 \Diff^+(K) \left| 
{\aligned
& \text{there is an element } 
\Phi \in \Diff^+ (S^4) \\ 
& \text{such that } 
\Phi|_K \text{ represents } \phi 
\endaligned}
\right. 
\right\} , 
\end{align*}
and define a quadratic form 
({\it the Rokhlin quadratic form\/}) 
$q_K\co H_1(K; {\Bbb Z}_2) \to {\Bbb Z}_2$: 
Let $P$ be a compact surface embedded in $S^4$, 
with its boundary contained in $K$, normal to $K$ along its 
boundary, and its interior is transverse to $K$. 
Let $P'$ be a surface transverse to $P$ obtained by sliding $P$ parallel to 
itself over $K$. 
Define $q_K([\partial P]) = \#( \text{int} P \cap (P' \cup K)) \text{
mod }  2$, where int means the interior. 
This is a well-defined quadratic form with respect to 
the ${\Bbb Z}_2$-homology intersection form $(,)_2$ on $K$, 
i.e. for each pair of elements $x$, $y$ of 
$H_1(K; {\Bbb Z}_2)$, $q_K(x+y) = q_K(x) + q_K(y) + (x,y)_2$. 
For the trivial $\Sigma_g$-knot $(S^4, \Sigma_g)$, 
let ${\cal SP}_g$ be the subgroup of ${\cal M}_g$ whose 
elements leave $q_{\Sigma_g}$ invariant. This group 
${\cal SP}_g$ is called the {\it spin mapping class group\/} 
\cite{Harer}. In the case when $g=1$, Montesinos showed: 
\begin{thm}\label{thm:Montesinos}{\rm\cite{Montesinos}}\qua
%
${\cal E}(S^4,\Sigma_1) = {\cal SP}_1$. 
\end{thm}
In this paper, we generalize this result to higher genus:
\begin{thm}\label{thm:main}
%
For any $g \geq 1$, ${\cal E}(S^4,\Sigma_g) = {\cal SP}_g$.   
\end{thm}
The group ${\cal E}(S^4,K)$ remains unknown for many non-trivial
$\Sigma_g$-knots $K$.  On the other hand, for some class of non-trivial
$\Sigma_1$-knots
$(S^4, K)$,  Iwase \cite{Iwase} and the author \cite{Hirose}
determined the groups 
${\cal E}(S^4, K)$. 

Finally, 
the author would like to express his gratitude to Professor Masahico 
Saito for his helpful comments, and to Professor Nariya Kawazumi 
for introducing him results of Johnson \cite{Johnson2}. 
This paper was written while the author 
stayed at Michigan State University as a visiting scholar 
sponsored by the Japanese Ministry of Education, Culture, Sports, 
Science and Technology. 
He is grateful to the Department of Mathematics, 
Michigan State University, for its hospitality. 
\section{Some elements of ${\cal E}(S^4, \Sigma_g)$}
For elements $a$, $b$ and $c$ of a group, 
we write 
$\inv{c} = c^{-1}$, and $a*b = ab \inv{a}$. 
Here, we introduce a standard form of the trivial $\Sigma_g$-knot 
$(S^4, \Sigma_g)$. 
We decompose $S^4 = D^4_+ \cup D^4_-$ and call 
$S^3 = D^4_+ \cap D^4_-$ {\it the equator $S^3$ \/}, and 
decompose $S^3 = D^3_+ \cup D^3_-$ and call 
$S^2 = D^3_+ \cap D^3_-$ {\it the equator $S^2$ \/}. 
\begin{figure}[ht!]
\centering
\includegraphics[height=4cm]{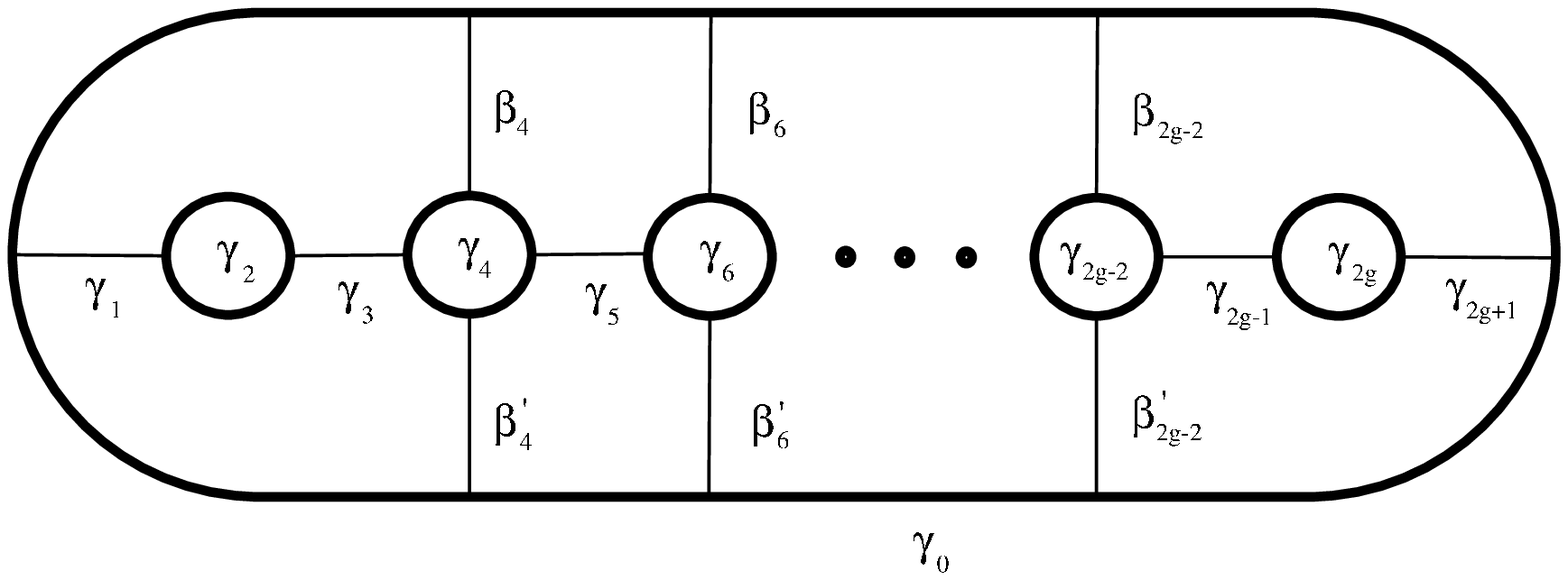}
\nocolon
\caption{}
\label{fig:planar}
\end{figure}
Let $P_g$ be a planar surface constructed from a 2-disk by removing 
$g$ copies of disjoint 2-disks. 
As indicated in Figure \ref{fig:planar}, we denote the boundary components of 
$P_g$ by $\gamma_0, \gamma_2, \ldots, \gamma_{2g}$, and denote 
some properly embedded arcs of $P_g$ by 
$\gamma_1, \gamma_3, \ldots, \gamma_{2g+1}$, 
$\beta_2, \beta_4, \ldots, \beta_{2g-2}$ and 
$\beta'_2, \beta'_4, \ldots, \beta'_{2g-2}$. 
We parametrize the regular neighborhood of the equator $S^2$ 
in the equator $S^3$ by $S^2 \times [-1,1]$, such that 
$S^2 \times \{0\}$ $=$ the equator $S^2$, 
$S^2 \times [-1,1] \cap D^3_+ = S^2 \times [0,1]$ and 
$S^2 \times [-1,1] \cap D^3_- = S^2 \times [-1,0]$. 
We put $P_g$ on the equator $S^2$. 
Then, $P_g \times[-1,1] \subset S^2 \times [-1,1]$ is a 3-dimensional 
handle body, so that, $(S^4, \partial (P_g \times[-1,1]))$ is 
the trivial $\Sigma_g$-knot. 
On $\partial (P_g \times[-1,1]) = \Sigma_g$, we define 
$c_{2i-1} = \partial (\gamma_{2i-1} \times [-1,1])$ $(1 \leq i \leq g+1)$, 
$b_{2j} = \partial (\beta_{2j}\times [-1,1])$, 
$b'_{2j} = \partial (\beta'_{2j}\times [-1,1])$ $(2 \leq j \leq g-1)$, and 
$c_{2k} = \gamma_{2k} \times \{ 0 \}$ $(1 \leq k \leq g)$. 
\begin{figure}[ht!]
\centering
\includegraphics[height=3cm]{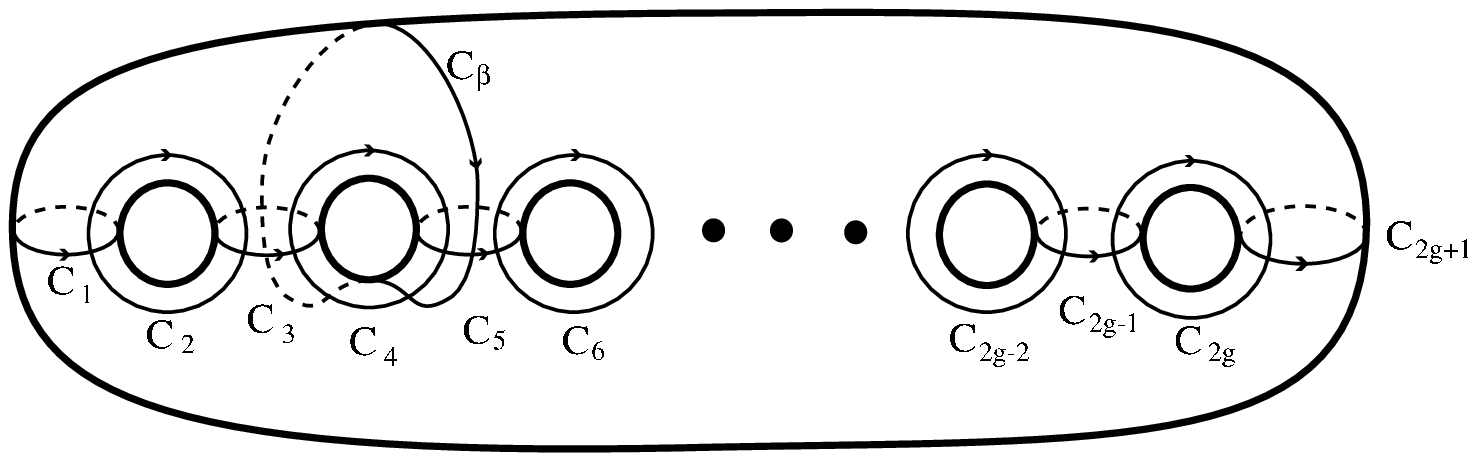}
\nocolon
\caption{}
\label{fig:circle}
\end{figure}
\begin{figure}[ht!]
\centering
\includegraphics[height=3.5cm]{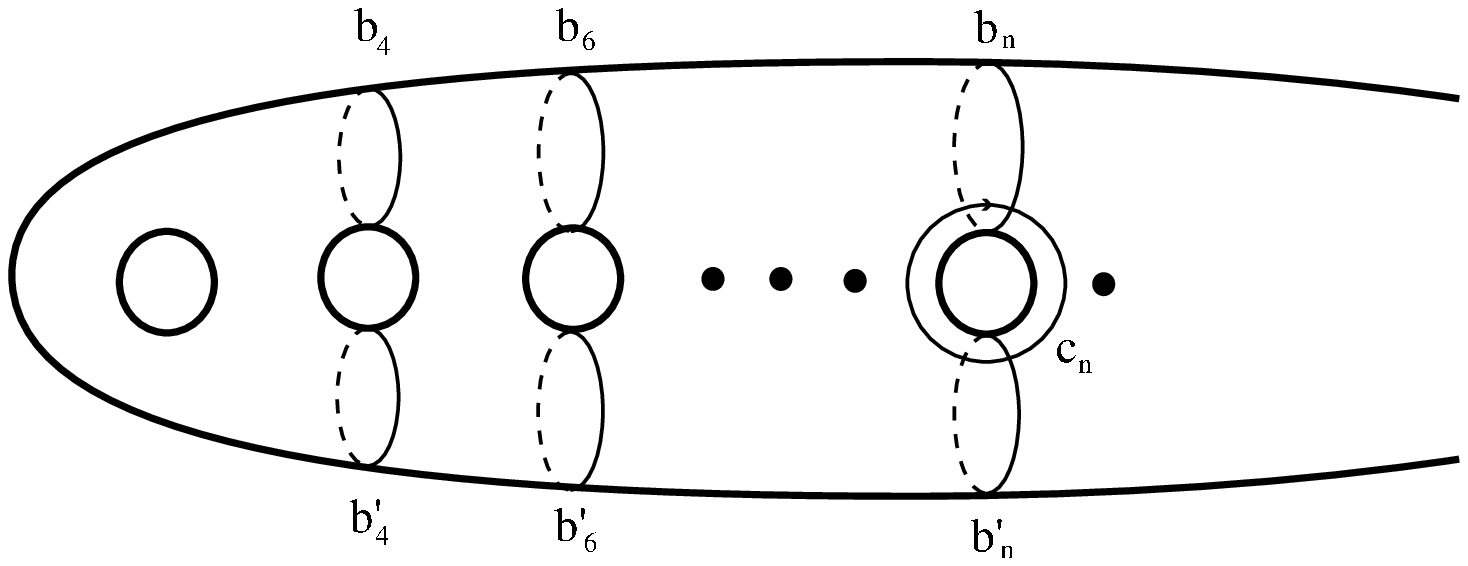}
\nocolon
\caption{}
\label{fig:b-circle}
\end{figure}
In Figures \ref{fig:circle} and \ref{fig:b-circle}, these circles are 
illustrated and some of them are oriented. 
For a simple closed curve $a$ on $\Sigma_g$, we denote 
the Dehn twist about $a$ by $T_a$. 
The order of composition of maps is the functional one: 
$T_b T_a$ means we apply $T_a$ first, then $T_b$. 
We define some elements of ${\cal M}_g$ as follows:
\begin{equation}
\begin{align*}
&C_i = T_{c_i},\  B_i = T_{b_i},\  B'_i = T_{b'_i}, \\
&X_i = C_{i+1} C_i \inv{C_{i+1}},\  
X^*_i = \inv{C_{i+1}} C_i C_{i+1} \ \ (1 \leq i \leq 2g),\\ 
&Y_{2j} = C_{2j} B_{2j} \inv{C_{2j}},\  
Y^*_{2j} = \inv{C_{2j}} B_{2j} C_{2j} \ \ (2 \leq j \leq g-1),\\ 
&D_i = C_i^2 \ \ (1 \leq i \leq 2g+1), \\ 
&DB_{2j} = B_{2j}^2 \ \ (2 \leq j \leq g-1), \\ 
&T=C_1 C_3 C_5,\  T_1 = C_1 C_3 B_4,\  
T_2 = B_4 C_5 C_7 \cdots C_{2g+1}.
\end{align*}
\end{equation}
When $g \geq 3$, 
the subgroup of ${\cal M}_g$ generated by 
$X_i$ $(1 \leq i \leq 2g)$, $Y_{2j}$ $(2 \leq j \leq g-1)$, 
$D_i$ $(1 \leq i \leq 2g+1)$, $DB_{2j}$ $(2 \leq j \leq g-1)$, 
$T_1$, and $T_2$ is denoted by $G_g$. 
It is clear that $X^*_i$ and $Y^*_{2j}$ are elements 
of $G_g$. 
When $g=2$, the subgroup of ${\cal M}_2$ generated by 
$X_i$ $(1 \leq i \leq 4)$, $D_j$ $(1 \leq j \leq 5)$, and 
$T$ is denoted by $G_2$. 
For two simple closed curves $l$ and $m$ on $\Sigma_g$, 
$l$ and $m$ are called {\it $G_g$-equivalent\/} (denote by 
$l \Ggeq m$) if there is an element $\phi$ of $G_g$ 
such that $\phi(l) = m$. 
\begin{figure}[ht!]
\centering
\includegraphics[height=2.5cm]{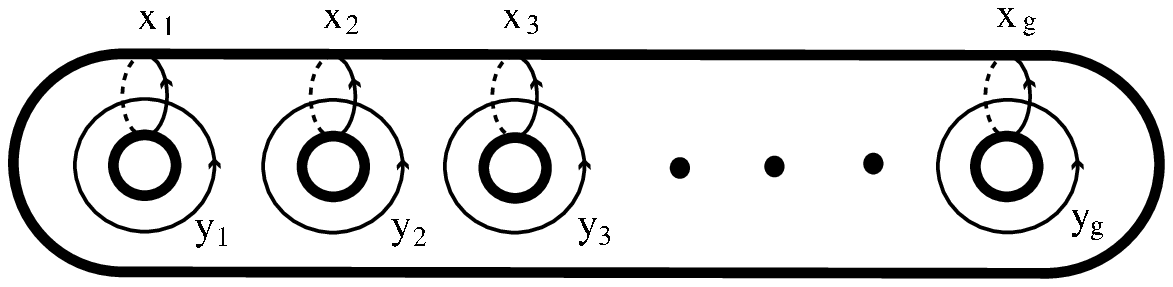}
\nocolon
\caption{}
\label{fig:basis}
\end{figure}
We set a basis of $H_1 (\Sigma_g; \Bbb Z)$ as in Figure \ref{fig:basis}, 
then for the quadratic form $q_{\Sigma_g}$ defined in \S 1, 
$q_{\Sigma_g}(x_i) = q_{\Sigma_g}(y_i) = 0$ ($1 \leq i \leq g$). 
By the definitions of $q_{\Sigma_g}$ and ${\cal SP}_g$, we have:
\begin{lem}\label{lem:EinSP} 
%
${\cal E}(S^4, \Sigma_g) \subset {\cal SP}_g$. 
\end{lem}
In this section, we show: 
\begin{lem}\label{lem:GinE} 
%
$G_g \subset {\cal E}(S^4, \Sigma_g)$. 
\end{lem}
As a straightforward corollary of these lemmas, we have: 
\begin{cor}\label{cor:GinSP}
%
%
$G_g \subset {\cal SP}_g$. 
\end{cor}
If $G_g \supset {\cal SP}_g$, then Theorem \ref{thm:main} is proved. 
We prove $G_g \supset {\cal SP}_g$ in the next section.
\begin{proof}[Proof of Lemma \ref{lem:GinE}] First we show that, if $g=2$,  
$T = C_1 C_3 C_5$ is an element of ${\cal E}(S^4, \Sigma_2)$. 
We parametrize the regular neighborhood of the equator $S^3$ in $S^4$ 
by $S^3 \times [-1,1]$, such that 
$S^3 \times \{ 0 \}$ $=$ the equator $S^3$, 
$S^3 \times [-1,1] \cap D^4_-$ $=$ $S^3 \times [-1,0]$, 
and  $S^3 \times [-1,1] \cap D^4_+$ $=$ $S^3 \times [0,1]$. 
\begin{figure}[ht!]
\centering
\includegraphics[height=9.5cm]{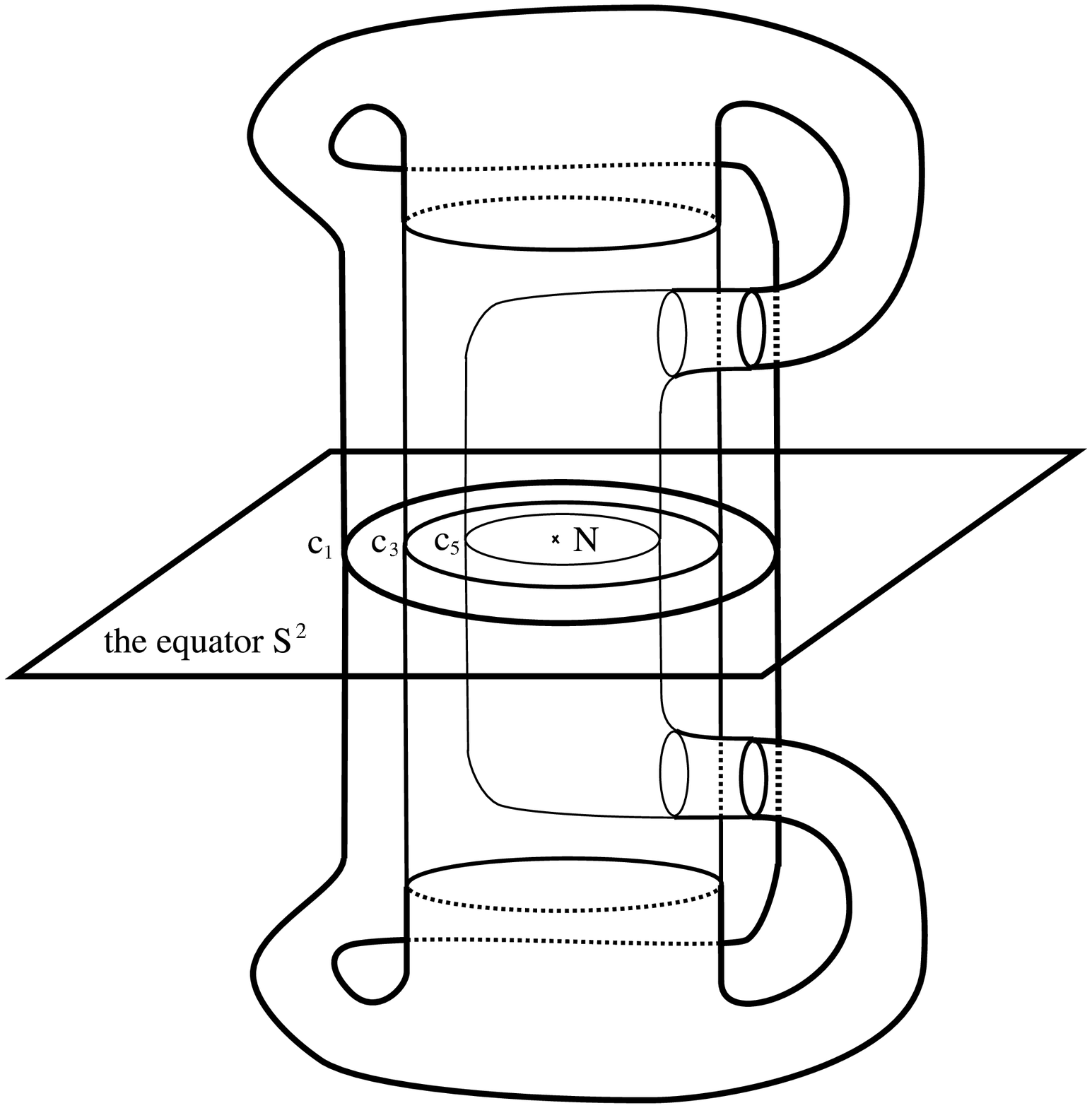}
\nocolon
\caption{}
\label{fig:concentric}
\end{figure}
We deform $\Sigma_2$ in $S^4$, in such a way that the surface obtained as a
result  of this deformation projects onto the equator $S^3$ 
as indicated in Figure \ref{fig:concentric}.  
In this figure, there are 6 intersecting circles.  
For each circle, we take two regular neighborhoods $N_1$ and $N_2$  in
$\Sigma_2$.  For $0< \epsilon < 1$, we put $N_1$ into 
$S^3 \times \{\frac{\epsilon}{2} \}$ and 
$N_2$ into $S^3 \times \{-\frac{\epsilon}{2} \}$. 
This deformation defines an orientation preserving diffeomorphism 
$\Psi_1$ of $S^4$. 
Let $r(\theta)\co S^2 \to S^2$ be the angle $\theta$ rotation whose axis 
passes through $N$. 
We define $R(\theta)\co S^3 \to S^3$ by 
\begin{equation}
\begin{align*}
R(\theta)&(x,t) = (r(t\theta)(x), t ) \quad 
\text{ on } S^2 \times [0,1] \\
R(\theta) & = id \quad \text{ on } D^3_- \\
R(\theta) & = \text{ the angle } \theta \text{ rotation } \quad 
\text{ on } D^3_+ - S^2 \times [0,1]. 
\end{align*}
\end{equation}
We define an orientation preserving diffeomorphism $\Psi_2$ of 
$S^4$ by 
\begin{equation}
\begin{align*}
\Psi_2 & (x,t) = (R(2 \pi)(x), t) \quad \text{ on } S^3 \times 
[-\epsilon, \epsilon], \\
\Psi_2 & (x,t) = \left( R(2 \pi \frac{1-t}{1-\epsilon})(x), t \right) \quad 
\text{ on } S^3 \times [\epsilon, 1], \\
\Psi_2 & (x,t) = \left( R(2 \pi \frac{t+1}{1-\epsilon})(x), t \right) \quad 
\text{ on } S^3 \times [-1,-\epsilon], \\
\Psi_2& = id \quad \text{ on } S^4 - S^3 \times [-1,1]. 
\end{align*}
\end{equation}
Then $\Psi_1^{-1} \Psi_2 \Psi_1 | _{\Sigma_2}$ $=$ $C_1 C_3 C_5$. 
In the same way as above, we can show for $g \geq 3$ that  
$T_1$ and  $T_2$ are elements of ${\cal E}(S^4, \Sigma_g)$. 

Next, for $g=3$, we show that $X_3 = C_4 C_3 \inv{C_4}$ and $D_3 = C_3^2$
are  elements of ${\cal E}(S^4, \Sigma_g)$. 
We review a theorem due to Montesinos \cite{Montesinos}. 
We can construct $S^4$ from $B^3 \times S^1$ and $S^2 \times D^2$ by 
attaching their boundary with the natural identification. 
Let $D^2 \times S^1$ be the solid torus trivially embedded in $B^3$. 
We regard $D^2 \times S^1 \times S^1$ $\subset$ $B^3 \times S^1$ 
$\subset$ $S^4$ as the regular neighborhood of a trivial $\Sigma_1$-knot. 
Let $E^4$ be the exterior of this trivial $\Sigma_1$-knot. 
The 3 simple closed curves $l=\partial D^2 \times * \times *$, 
$r=* \times S^1 \times *$, $s=* \times * \times S^1$ on 
$\partial E^4$ represent a basis of $H_1(\partial E^4; {\Bbb Z})$. 
Montesinos showed: 
\begin{thm}\label{Montesinos}{\rm\cite[Theorem 5.3]{Montesinos}}\qua
%
%
Let $g \co \partial E^4 \to \partial E^4$ be a diffeomorphism which induces 
an automorphism on $H_1(\partial E^4; {\Bbb Z})$, 
$$
g_*(l,r,s) = (l,r,s)
\begin{pmatrix}
m & a & b \\
n & \alpha & \gamma \\
p & \beta & \delta 
\end{pmatrix}.
$$
There is a diffeomorphism $G \co E^4 \to E^4$ such that $G|_{\partial E^4} = g$ 
if and only if $a=b=0$ and $\alpha + \beta + \gamma + \delta$ is even. 
\end{thm}
Let $p$ be a point on $* \times S^1 \times S^1$ disjoint from $r \cup s$, 
$N(p)$ be a regular neighborhood of $p$ in the equator $S^3$, then 
$N$ $=$ $* \times S^1 \times S^1 - N(p)$ in a regular neighborhood of 
$r \cup s$. 
\begin{figure}[ht!]
\centering
\includegraphics[height=3.8cm]{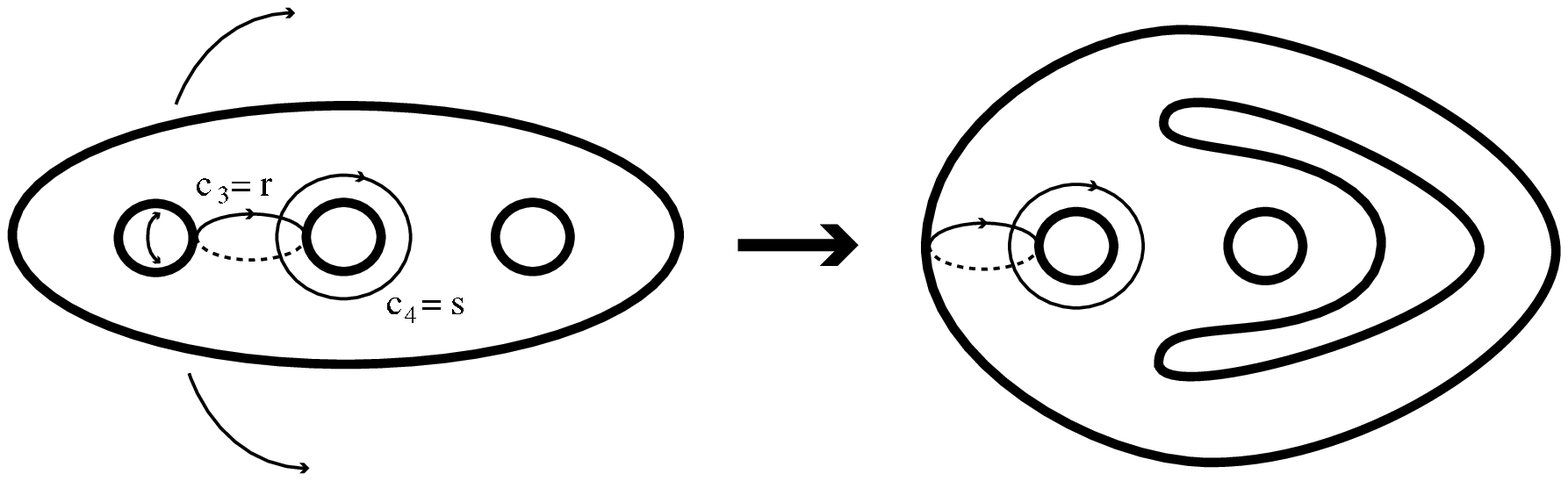}
\nocolon
\caption{}
\label{fig:deformS3}
\end{figure}
Figure \ref{fig:deformS3} illustrates deformation of $\Sigma_g$ into 
$D^2 \times S^1 \times S^1$. 
We bring $c_3$ and $c_4$ to $r$ and $s$ and deform as is indicated 
by arrows. 
Then, we can deform $\Sigma_3$ in such a way that a regular neighborhood
$N'$ of 
$c_3 \cup c_4$ coincides with $N$ and 
$\Sigma_3 - N'$ $\subset$ $N(p)$. 
Let diffeomorphisms $f_1$, $f_2$ over $D^2 \times S^1 \times S^1$ 
be defined by 
$f_1 = id_{D^2} \times 
\begin{pmatrix}
1 & 2 \\
0 & 1 
\end{pmatrix}
$, 
$f_2 = id_{D^2} \times 
\begin{pmatrix}
2 & 1 \\
-1 & 0 
\end{pmatrix}
$ 
(where we present diffeomorphisms on $* \times S^1 \times S^1$ by its action 
on the basis $\{r, s\}$ of $H_1(*\times S^1 \times S^1;{\Bbb Z})$ and 
$r$ and $s$ are oriented as in Figure \ref{fig:deformS3}), 
then $f_1 |_{\Sigma_2}$ $=$ $C_3^2= D_3$, 
$f_2|_{\Sigma_2}$ $=$ $C_4 C_3 \inv{C_4} = X_3$. 
Since the actions of these homeomorphisms on $H_1(\partial E^4; {\Bbb Z})$
are described by  
$$
(f_1 | \partial E^4)_*(l,r,s) = (l,r,s) 
\begin{pmatrix}
1 & 0 & 0 \\
0 & 1 & 2 \\
0 & 0 & 1 
\end{pmatrix},
$$
$$
(f_2 | \partial E^4)_*(l,r,s) = (l,r,s) 
\begin{pmatrix}
1 & 0 & 0 \\
0 & 2 & 1 \\
0 & -1 & 0 
\end{pmatrix},
$$
there are diffeomorphisms $F_1$ and $F_2$ such that 
$F_1 | _{D^2 \times S^1 \times S^1} = f_1$, 
$F_2 | _{D^2 \times S^1 \times S^1}$  $=$ $f_2$. 
These diffeomorphisms $F_1$, $F_2$ are extensions of 
$f_1$, $f_2$ respectively. 
By the same method as above, we can show that other $X_i$, $Y_{2j}$, $D_i$, 
and $DB_{2j}$ are elements of ${\cal E}(S^4, \Sigma_g)$ 
for any $g \geq 2$. 
\end{proof}
\section{A finite set of generators for the spin mapping class group}
In Corollary \ref{cor:GinSP}, we showed that $G_g \subset {\cal SP}_g$. 
In this section, we show that $G_g = {\cal SP}_g$. That is to say, 
we show: 
\begin{thm}\label{thm:generator-Gg}
%
If $g=2$, ${\cal SP}_2$ is generated by 
$C_{i+1} C_i \inv{C_{i+1}}$ ($1 \leq i \leq 4$), $C_j^2$ ($1 \leq j \leq 5$), 
and $C_1 C_3 C_5$. 
If $g \geq 3$, $ {\cal SP}_g$ is generated by 
$C_{i+1} C_i \inv{C_{i+1}}$ ($1 \leq i \leq 2g$), 
$C_{2j} B_{2j} \inv{C_{2j}}$ ($2 \leq j \leq g-1$), 
$C_k^2$ ($1 \leq k \leq 2g+1$), 
$B_l^2$ ($1 \leq l \leq g-1$), $C_1 C_3 B_4$ and $B_4 C_5 C_7 \cdots C_{2g+1}$. \end{thm}
When $g=2$, we use Reidemeister--Schreier's method to show this. 
On the other hand, when $g \geq 3$, we use other methods. 
We start from the case when $g \geq 3$. 
\subsection{The hyperelliptic mapping class group}
Let ${\cal H}_g$ be the subgroup of the mapping class group 
${\cal M}_g$ generated by $C_1, C_2,$ $\ldots,$ $C_{2g+1}$. 
This group is called {\it the hyperelliptic mapping class group\/}. 
In this group (and also in ${\cal M}_g$), 
$C_i$'s satisfy the following equations:
\begin{equation}
\begin{align*}
C_i C_{i+1} C_i &= C_{i+1}  C_i C_{i+1}, \; (1 \leq i \leq 2g) \qquad \qquad \\
C_i C_j &= C_j C_i, \; (|i-j| \geq 2). \qquad \qquad
\end{align*}
\end{equation}
These equations are called {\it braid equation\/}. 
In this paper, we use these relations frequently. 
In this section, we show the following lemma for ${\cal H}_g$. 
\begin{lem}\label{lem:Braid-Conj}
%
%
For any $i = 1,2,\ldots,2g+1$, and any element $W$ of ${\cal H}_g$, 
$W C_i C_i \inv{W}$ is an element of $G_g$. 
\end{lem}
\begin{proof}
We call $C_i$ {\it a positive letter\/} and 
$\inv{C_i}$ {\it a negative letter\/}. 
A sequence of positive letters is called {\it a positive word}. 
If indices of two letters $C_i$, $C_j$ satisfy $|i-j|=1$, 
then we say $C_i$ is {\it adjacent\/} to $C_j$. 
If there is a negative letter $\inv{B}$ in a sequence of letters 
$W$, which presents an element of ${\cal H}_g$, we replace 
$\inv{B}$ by a sequence of letters $\inv{B} \inv{B} \cdot B$. 
This shows that every element of ${\cal H}_g$ is 
represented by a sequence of positive letters and 
$\inv{C_j}\inv{C_j}$'s $(1\leq j \leq 2g+1)$. 
If there is a sequence of letters $XX$ ( $X=C_i$ or $\inv{C_i}$) 
in $W$, say $W$ $=$ $W_1 XX W_2$, then we rewrite, 
\begin{equation}
\begin{align*}
W C_i C_i \inv{W} 
&= W_1 XX W_2 C_i C_i \inv{W_2} \inv{X} \inv{X} 
\inv{W_1} \\
&= W_1 XX \inv{W_1} W_1 W_2 C_i C_i \inv{W_2}
\inv{W_1} W_1 \inv{X}\inv{X} \inv{W_1}.
\end{align*}
\end{equation}
Therefore, the following claim shows this lemma:
\par\medskip
{\bf Claim}\qua{\sl For any positive word $W$ without 
$C_j C_j (1 \leq j \leq 2g+1)$, 
$W C_i C_i \inv{W}$ is an element of $G_g$. \/} 
\par\medskip
If the word length of $W$ is 0, the above claim is trivial. 
We assume that the word length of $W$ is at least 1, and we show this 
claim by the induction on the word length. 
If the right most letter $L$ of $W$ is not adjacent to $A_i$, 
and say $W = W'L$, then 
$$
WC_i C_i \inv{W} = W' L C_i C_i \inv{L} \inv{W'} 
=W'C_i L \inv{L} C_i \inv{W'} 
=W'C_i C_i \inv{W'}.
$$
By the induction hypothesis, $W C_i C_i \inv{W}$ is 
an element of $G_g$. 
Therefore, from here to the end of this proof, we assume 
that the right most letter of $W$ is adjacent to $C_i$. 
Let $l$ be the word length of $W$, and  
$W$ $= x_l x_{l-1} \ldots x_2 x_1$. 
The letter $x_i$ of $W$ is called {\it a jump\/}, if 
$x_{i-1}$ and $x_i$ are not adjacent. 
The letter $x_j$ of $W$ is called {\it a turn\/}, if 
$x_j$ and $x_{j-1}$ are not jumps and $x_j = x_{j-2}$. 
Considering jumps and turns, we need to show this 
claim for the following three cases.

\medskip
{\bf Case 1}\qua{\it When there is not any jump or any turn\/}: 
Since $x_l$ and $x_{l-1}$ are adjacent, 
$x_l x_{l-1} \inv{x_l}$ is an element of $G_g$. 
We rewrite, 
$$
WC_i C_i \inv{W} 
= x_l x_{l-1} \inv{x_l} 
\cdot x_l x_{l-2} x_{l-3} \cdots x_1 C_i C_i 
\inv{x_1} \cdots \inv{x_{l-3}} \inv{x_{l-2}} 
\inv{x_l} 
\cdot x_l \inv{x_{l-1}} \inv{x_l}. 
$$
By the induction hypothesis, $W C_i C_i \inv{W}$ is 
an element of $G_g$. 

\medskip
{\bf Case 2}\qua{\it When there are jumps, but there is not any turn\/}: 
We show in the induction on the number of jumps in $W$. 
Let $x_j$ be the right most jump in $W$. 
First we consider the case when $j = 2$, say 
$W = W' x_2 x_1$. If $x_2$ is not adjacent to $C_i$, we 
rewrite,
\begin{equation}
\begin{align*}
W C_i C_i \inv{W} &= W' x_2 x_1 C_i C_i \inv{x_1} 
\inv{x_2} \inv{W'} \qquad \qquad \qquad \\
&= W' x_1 x_2 C_i C_i \inv{x_2} \inv{x_1} \inv{W'} \\
&= W' x_1 C_i x_2 \inv{x_2} C_i \inv{x_1} \inv{W'} \\
&= W' x_1 C_i C_i \inv{x_i} \inv{W'}. 
\end{align*}
\end{equation}
By the induction hypothesis on the word length of $W$, 
$W C_i C_i \inv{W}$ is an element of $G_g$. 
If $x_2$ is adjacent to $C_i$, we rewrite, 
\begin{equation}
\begin{align*}
W C_i C_i \inv{W} &= W' x_2 x_1 C_i C_i \inv{x_1} \inv{x_2} \inv{W'} \\ 
&= W' x_2 \inv{C_i} x_1 x_1 C_i \inv{x_2} \inv{W'} \\
&= W' x_2 \inv{C_i} \inv{C_i} \cdot C_i x_1 x_1 \inv{C_i} 
\cdot C_i C_i \inv{x_2} \inv{W'} \\
&= W' x_2 \inv{C_i} \inv{C_i} \inv{x_2} \inv{W'} 
\cdot W' x_2 C_i x_1 x_1 \inv{C_i} \inv{x_2} \inv{W'} 
\cdot W' x_2 C_i C_i \inv{x_2} \inv{W'}.  
\end{align*}
\end{equation}
By the induction hypothesis on the word length of $W$, 
the first and third terms are elements of $G_g$. 
By the induction hypothesis on the number of jumps in $W$, 
the second term is an element of $G_g$. 
Therefore, $W C_i C_i \inv{W}$ is an element of $G_g$. 
Next, we consider on the case when $j$ is at least 3. 
If $x_j$ is not adjacent to $x_{j-1}, \ldots, x_1$ then, 
$$
W = \ldots x_j x_{j-1} \ldots x_1 = \ldots x_{j-1} \ldots x_1 x_j. 
$$
Therefore, it comes down to the case $j=2$. 
If there are some letters adjacent to $x_j$ in 
$\{ x_{j-1}, \cdots, x_1 \}$, let $x_i$ be the left most element 
among them. 
By the definition of jumps, $j > i+1$, and by the definition of 
$x_i$, $x_j=x_{i-1}$. Therefore, 
\begin{equation}
\begin{align*}
W &= \cdots x_j \cdots x_{i+1} x_i x_{i-1} \cdots x_1 \qquad \qquad \qquad\\
&= \cdots x_{i+1} x_j x_i x_{i-1} \cdots x_1 \\
&= \cdots x_{i+1} x_{i-1} x_i x_{i-1} \cdots x_1 \\
&= \cdots x_{i+1} x_i x_{i-1} x_i \cdots x_1. 
\end{align*}
\end{equation}
Since there is not any jump or any turn in the sequence 
$x_i x_{i-1} \cdots x_1$, $x_i$ commutes with 
$x_{i-2}, \ldots, x_1$. Therefore, 
$W = \cdots x_1 x_i$ and it comes down to the case $j=2$. 

\medskip
{\bf Case 3}\qua  {\it When there are turns in $W$\/}: 
Let $x_t$ be the right most turn in $W$. By the definition of 
turn, $t$ is at least 3. 
By applying the argument for Case 2 to $x_{t-1} x_{t-2} \cdots x_1$,  
we assume that there is no turn and no jump in 
$x_{t-1} x_{t-2} \cdots x_1$. 
Since we assume that $x_1$ is adjacent to $C_i$, there may be  
a case when $x_2 = C_i$. 
In that case, we rewrite, 
\begin{equation}
\begin{align*}
W C_i C_i \inv{W} &= \cdots x_3 x_2 x_1 C_i C_i 
\inv{x_1} \inv{x_2}\inv{x_3}\cdots \\
&= \cdots x_3 C_i x_1 C_i C_i \inv{x_1} \inv{C_i} \inv{x_3} \cdots \qquad \qquad \qquad\\
&= \cdots x_3 x_1 C_i x_1 \inv{x_1} \inv{C_i} x_1 \inv{x_3} \cdots \\
&= \cdots x_3 x_1 x_1 \inv{x_3} \cdots. 
\end{align*}
\end{equation}
By the induction hypothesis on the word length of $W$, 
$W C_i C_i \inv{W}$ is an element of $G_g$. 
If $x_2 \not= C_i$, then  
$x_{t-1}, x_{t-2}, \cdots, x_2$ are not adjacent to $C_i$. 
We rewrite, 
\begin{equation}
\begin{align*}
W &= \cdots x_t x_{t-1} x_{t-2} x_{t-3} \cdots x_1 \qquad \qquad \qquad\\
&= \cdots x_{t-2} x_{t-1} x_{t-2} x_{t-3} \cdots x_1 \\
&= \cdots x_{t-1} x_{t-2} x_{t-1} x_{t-3} \cdots x_1. 
\end{align*}
\end{equation}
Since we assume that there is no jump and no turn in 
$x_{t-1} x_{t-2} \cdots x_1$, $x_{t-1}$ is not adjacent to 
$x_{t-3}, \ldots, x_1$. 
Therefore, $W = \cdots x_{t-1} x_{t-2} x_{t-3} \cdots x_1 x_{t-1}$. 
With remarking that $x_{t-1}$ is not adjacent to $C_i$, we rewrite, 
\begin{equation}
\begin{align*}
W C_i C_i \inv{W} 
&= \cdots x_{t-1} x_{t-2} x_{t-3} \cdots x_1 x_{t-1} C_i C_i 
\inv{x_{t-1}} \inv{x_1} \cdots \inv{x_{t-3}} \inv{x_{t-2}} \inv{x_{t-1}} 
\cdots \\
&= \cdots x_{t-1} x_{t-2} x_{t-3} \cdots x_1 C_i x_{t-1} \inv{x_{t-1}} 
C_i \inv{x_1} \cdots \inv{x_{t-3}} \inv{x_{t-2}} \inv{x_{t-1}} 
\cdots \\   
&= \cdots x_{t-1} x_{t-2} x_{t-3} \cdots x_1 C_i C_i \inv{x_1} \cdots
\inv{x_{t-3}} \inv{x_{t-2}} \inv{x_{t-1}} \cdots. 
\end{align*}
\end{equation}
By the induction hypothesis on the word length of $W$, 
$W C_i C_i \inv{W}$ is an element of $G_g$. 
\end{proof}
\subsection{The Torelli group ${\cal I}_g$}  
\begin{figure}[ht!]
\centering
\includegraphics[height=1.4cm]{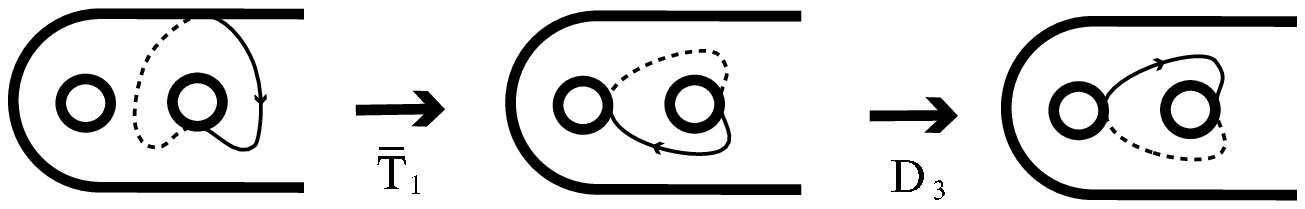}
\nocolon
\caption{}
\label{fig:BtoS}
\end{figure}

In this subsection, we assume $g \geq 3$.  
There is a natural surjection $\Phi\co {\cal M}_g \to \Symp (2g, {\Bbb Z})$ 
defined by the action of ${\cal M}_g$ on the group 
$H_1(\Sigma_g; {\Bbb Z})$. We denote the kernel of $\Phi$ by 
${\cal I}_g$ and call this {\it the Torelli group\/}. 
In this subsection, we prove the following lemma:
\begin{lem}\label{lem:Torelli}
%
The Torelli group ${\cal I}_g$ is a subgroup of 
$G_g$. 
\end{lem}
\begin{figure}[ht!]
\centering
\includegraphics[height=1.5cm]{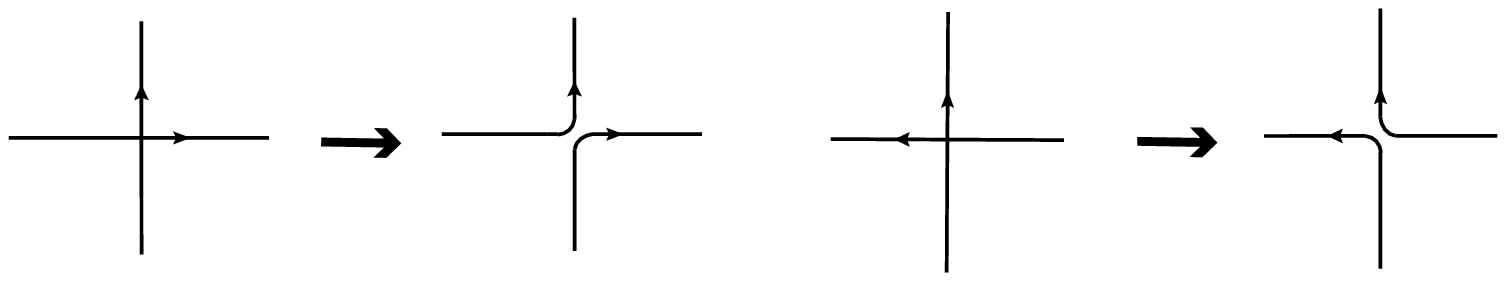}
\nocolon
\caption{}
\label{fig:smoothing}
\end{figure}
Johnson \cite{Johnson1} showed that, when $g$ is larger than or equal to 
$3$, ${\cal I}_g$ is finitely generated. We review his result. 
We orient and call simple closed curves as indicated in 
Figure \ref{fig:circle}, and call $(c_1, c_2,\ldots, c_{2g+1})$ and 
$(c_{\beta}, c_5, \ldots, c_{2g+1})$ as {\it chains\/}. 
For oriented simple closed curves $d$ and $e$ which mutually intersect 
in one point, we construct an oriented simple closed curve $d+e$ 
from $d \cup e$ as follows: 
choose a disk neighborhood of the intersection point and in it 
make a replacement as indicated in Figure \ref{fig:smoothing}. 
For a consecutive subset $\{ c_i, c_{i+1}, \ldots, c_j \}$ 
of a chain, let $c_i + \cdots + c_j$ be the oriented simple closed 
curve constructed by repeated applications of the above operations. 
Let $(i_1, \ldots, i_{r+1})$ be a subsequence of 
$(1, 2,\ldots, 2g+1)$ (Resp. $(\beta, 5, \ldots, 2g+1)$). 
We construct the union of circles 
${\cal C}$ $=$ $c_{i_1}+\cdots+c_{i_2-1} \cup c_{i_2}+\cdots+c_{i_3-1} 
\cup \cdots \cup c_{i_r}+\cdots+c_{i_{r+1}-1}$. 
If $r$ is odd, the regular neighborhood of ${\cal C}$ is an oriented 
compact surface with $2$ boundary components. 
Let $\phi$ be the element of ${\cal M}_g$ defined as 
the composition of 
the positive Dehn twist along the boundary curve to the left of ${\cal C}$
and  the negative Dehn twist along the boundary curve to the right of 
${\cal C}$.  Then, $\phi$ is an element of ${\cal I}_g$. We denote 
$\phi$ by $[i_1, \ldots, i_{r+1}]$, and call this 
{\it the odd subchain map\/} of $(c_1, c_2,\ldots, c_{2g+1})$ 
(Resp. $(c_{\beta}, c_5, \ldots, c_{2g+1})$). 
Johnson \cite{Johnson1} showed the following theorem: 
\begin{thm}\label{thm:Johnson-gen}{\rm\cite[Main
Theorem]{Johnson1}}\qua 
%
For $g \geq 3$, the odd subchain maps of the two chains 
$(c_1, c_2,\ldots, c_{2g+1})$ and $(c_{\beta}, c_5, \ldots, c_{2g+1})$
generate ${\cal I}_g$. 
\end{thm}
We use the following results by Johnson \cite{Johnson1}.  
\begin{lem}\label{lem:Johnson-action}{\rm\cite{Johnson1}}\qua 
%
{\rm(a)}\qua $C_j$ commutes with $[i_1, i_2, \cdots]$ if and only if $j$ and 
$j+1$ are either both contained in or are disjoint from the $i$'s. 
\newline
{\rm(b)}\qua If $i \not= j+1 $, then $\inv{C_j} * [ \cdots, j, i, \cdots ]$ 
$=$ $[ \cdots, j+1, i, \cdots ]$, and 
$C_j * [ \cdots, j, i, \cdots ]$ $=$ $[ \cdots, j, i, \cdots ] 
[ \cdots, j+1, i, \cdots]^{-1} [ \cdots, j, i, \cdots ]$. 
\newline
{\rm(c)}\qua If $k \not= j$, then $C_j*[\cdots, k, j+1, \cdots]$ $=$ 
$[\cdots, k, j, \cdots]$, and 
$\inv{C_j}*[\cdots, k, j+1, \cdots]$ $=$ 
$[\cdots, k, j+1, \cdots][\cdots, k, j, \cdots]^{-1}
[\cdots, k, j+1,\cdots]$. 
\newline
{\rm(d)}\qua $[1,2,3,4] [1,2,5,6,\ldots, 2n] B_4*[3,4,5,\ldots,2n]$ $=$ 
$[5,6,\ldots,2n][1,2,3,4,\ldots,
$ $
2n]$, where $3 \leq n \leq g$.      
\end{lem}
First we show that some odd subchain maps are elements of $G_g$. 
\begin{lem}\label{lem:special-elements}
%
$[1,2,3,4]$, $[1,3,5,7,\ldots, 2i+1, \ldots, 2n-1]$ 
$(n \text{ is even, and } 4\leq n \leq g+1)$, and 
$[1,2,4,6,\ldots,2i,\ldots,2n-2]$ 
$(n \text{ is even, and } 4\leq n \leq g+2)$ are elements of 
$G_g$. 
\end{lem}
\begin{proof}
In this proof, for a sequence $\{ f_i \}$ of elements of ${\cal M}_g$, 
we write, 
$$
\prod_{i=n}^m f_i = 
\begin{cases} 
f_n f_{n+1} \cdots f_m, & \qquad n \leq m, \\
f_n f_{n-1} \cdots f_m, & \qquad n \geq m. 
\end{cases}
$$
(1) {\it $[1, 2, 3, 4]$ is an element of $G_g$\/}: 
$[1,2,3,4]$ is equal to $B_4 \inv{B'_4}$. 
\linebreak
Since $C_4 C_3 C_2 C_1 C_1 C_2 C_3 C_4 (b_4)$ $=$ $b'_4$, 
\begin{equation}
\begin{align*}
[1,2,3,4] &= B_4 C_4 C_3 C_2 C_1 C_1 C_2 C_3 C_4 \inv{B_4} 
 \inv{C_4} \inv{C_3} \inv{C_2} \inv{C_1} \inv{C_1} \inv{C_2} 
\inv{C_3} \inv{C_4} \\
&= B_4 C_4 \inv{B_4} \cdot C_3 C_2 \inv{C_3} \cdot C_1 C_1 
\cdot C_3 C_2 \inv{C_3} \cdot C_3 C_3 \cdot B_4 C_4 \inv{B_4} 
\cdot \inv{C_4} \inv{C_3} C_4 \cdot \\ 
&\; \; \cdot \inv{C_2} \inv{C_1} C_2 \cdot 
\inv{C_2} \inv{C_1} C_2 \cdot \inv{C_2} \inv{C_2} \cdot
\inv{C_4} \inv{C_3} C_4 \cdot \inv{C_4} \inv{C_4}. 
\end{align*}
\end{equation}
Therefore, $[1,2,3,4]$ is an element of $G_g$. 
\newline
(2) {\it $[1,3,5,7,\ldots, 2i+1, \ldots, 2n-1]$ 
$(n \text{ is even, and } 4\leq n \leq g+1)$ are elements of $G_g$\/}: 
By (b) of Lemma \ref{lem:Johnson-action}, 
$$
[1,3,5,7,\ldots, 2i+1, \ldots, 2n-1] = 
(\prod_{k=n-1}^1 \prod_{i=2k}^{n+k-1} \inv{C_i}) 
* [1,2,3,4,\ldots,n]. 
$$
Since $[1,2,3,4,\ldots,n] = B_n \inv{B'_n}$, and 
$b'_n = \prod_{i=n}^2 C_i \cdot C_1 C_1 \cdot \prod_{i=2}^n C_i(b_n)$, 
\begin{equation}
\begin{align*}
[1,2,3,4,\cdots,n] &= B_n \prod_{i=n}^2 C_i \cdot C_1 C_1 \cdot
\prod_{i=2}^n C_i \cdot \inv{B_n} \cdot \prod_{i=n}^2 \inv{C_i}
\cdot \inv{C_1} \inv{C_1} \cdot \prod_{i=2}^n \inv{C_i} \\
&=\prod_{k=2}^n \{(B_n \prod_{i=n}^k C_i)*(C_{k-1} C_{k-1})\} 
\cdot B_n*(C_n C_n) \cdot \\
&\qquad \cdot \Prod{k=2}{n}{\{(\Prod{i=n}{k}{\inv{C_i}})*
(\inv{C_{k-1}} \inv{C_{k-1}}) \}} 
\cdot \inv{C_n} \inv{C_n}. 
\end{align*}
\end{equation}
Therefore, 
\begin{equation}
\begin{align*}
[1,3,5,7, \ldots, 2n-1] &= 
\Prod{k=2}{n}{ \{ (\Prod{l=n-1}{1}{\Prod{i=2l}{n+l-1}{\inv{C_i}}}
\cdot B_n \cdot \Prod{i=n}{k}{C_i}) * (C_{k-1} C_{k-1}) \} } \cdot \\
& \cdot (\Prod{l=n-1}{1}{\Prod{i=2l}{n+l-1}{\inv{C_i}}} \cdot B_n )
*(C_n C_n)\cdot \\
& \cdot \Prod{k=2}{n}{ \{ (\Prod{l=n-1}{1}{\Prod{i=2l}{n+l-1}{\inv{C_i}}}
\cdot \Prod{i=n}{k}{\inv{C_i}}) * (\inv{C_{k-1}} \inv{C_{k-1}}) \} }
\cdot \\
& \cdot (\Prod{l=n-1}{1}{\Prod{i=2l}{n+l-1}{\inv{C_i}}})
*(\inv{C_n} \inv{C_n}). 
\end{align*}
\end{equation} 
By Lemma \ref{lem:Braid-Conj}, $\Prod{k=2}{n}{ \{
(\Prod{l=n-1}{1}{\Prod{i=2l}{n+l-1}{\inv{C_i}}}
\cdot \Prod{i=n}{k}{\inv{C_i}}) * (\inv{C_{k-1}} \inv{C_{k-1}}) \} }$ and 
$(\Prod{l=n-1}{1}{\Prod{i=2l}{n+l-1}{\inv{C_i}}})*(\inv{C_n} \inv{C_n})$
are elements of $G_g$. 
By braid relations for ${\cal M}_g$, (in the following equations $j \leq
n-1$) 
\begin{equation}
\begin{align*}
(C_{j-1} \cdot &\Prod{i=n}{j}{C_i})*(C_{j-1} C_{j-1}) =
C_{j-1} \Prod{i=n}{j+1}{C_i} \cdot C_j C_{j-1} C_{j-1} \inv{C_j} \cdot 
\Prod{i=j+1}{n}{\inv{C_i}} \inv{C_{j-1}} \\
&= \Prod{i=n}{j+1}{C_i} \cdot C_{j-1} C_j C_{j-1} C_{j-1} \inv{C_j}
\inv{C_{j-1}} \cdot \Prod{i=j+1}{n}{\inv{C_i}}  \\
&= \Prod{i=n}{j+1}{C_i} \cdot C_j C_{j-1} C_j \inv{C_j} \inv{C_{j-1}}
C_j \cdot \Prod{i=j+1}{n}{\inv{C_i}}  
= (\Prod{i=n}{j+1}{C_i})*(C_j C_j),
\end{align*}
\end{equation}
\begin{equation}
\begin{align*}
(C_{n-1} C_n)*(C_{n-1} C_{n-1}) 
&= C_{n-1} C_n C_{n-1} C_{n-1} \inv{C_n} \inv{C_{n-1}} \qquad \qquad \qquad\\
&= C_n C_{n-1} C_n \inv{C_n} \inv{C_{n-1}} C_n 
= C_n C_n. 
\end{align*}
\end{equation}
By the above equation and the fact that $B_n$ commutes with 
$C_j$ ($1 \leq j \leq n-1)$, 
$$
(B_n \cdot \Prod{i=n}{k}{C_i})*(C_{k-1} C_{k-1}) = 
(\Prod{j=k-2}{1}{C_j} \cdot B_n \cdot \Prod{i=n}{2}{C_i}) 
*(C_1 C_1) 
\text{ where } 3 \leq k \leq n, 
$$
$$
B_n * (C_n C_n) = (\Prod{j=n-1}{1}{C_j} \cdot B_n 
\Prod{i=n}{2}{C_i}) * (C_1 C_1). 
$$
Since, for $3 \leq k \leq n+1$, 
$$
\Prod{l=n-1}{1}{\Prod{i=2l}{n+l-1}{\inv{C_i}}} \cdot 
\Prod{j=k-2}{1}{C_j} = 
\Prod{j=k-2}{1}{(\inv{C_{2j}} C_{2j-1} C_{2j})} \cdot
\Prod{l=n-1}{1}{\Prod{i=2l}{n+l-1}{\inv{C_i}}},    
$$
we obtain, 
\begin{equation}
\begin{align*}
(\Prod{l=n-1}{1}{\Prod{i=2l}{n+l-1}{\inv{C_i}}} \cdot 
&\Prod{j=k-2}{1}{C_j} \cdot B_n \cdot \Prod{i=n}{2}{C_i} )
* (C_1 C_1)  \\
&= 
(\Prod{j=k-2}{1}{(\inv{C_{2j}} C_{2j-1} C_{2j})} \cdot 
\Prod{l=n-1}{1}{\Prod{i=2l}{n+l-1}{\inv{C_i}}} \cdot 
B_n \cdot \Prod{i=n}{2}{C_i}) * (C_1 C_1). 
\end{align*}
\end{equation}
Therefore, for showing that $[1,3,5,7,\ldots,2n-1]$ is an  
element of $G_g$, it suffices to show that 
$(\Prod{l=n-1}{1}{\Prod{i=2l}{n+l-1}{\inv{C_i}}} \cdot
B_n \cdot \Prod{i=n}{2}{C_i})*(C_1 C_1)$ 
is an element of $G_g$.
\begin{figure}[ht!]
\centering
\includegraphics[height=3.2cm]{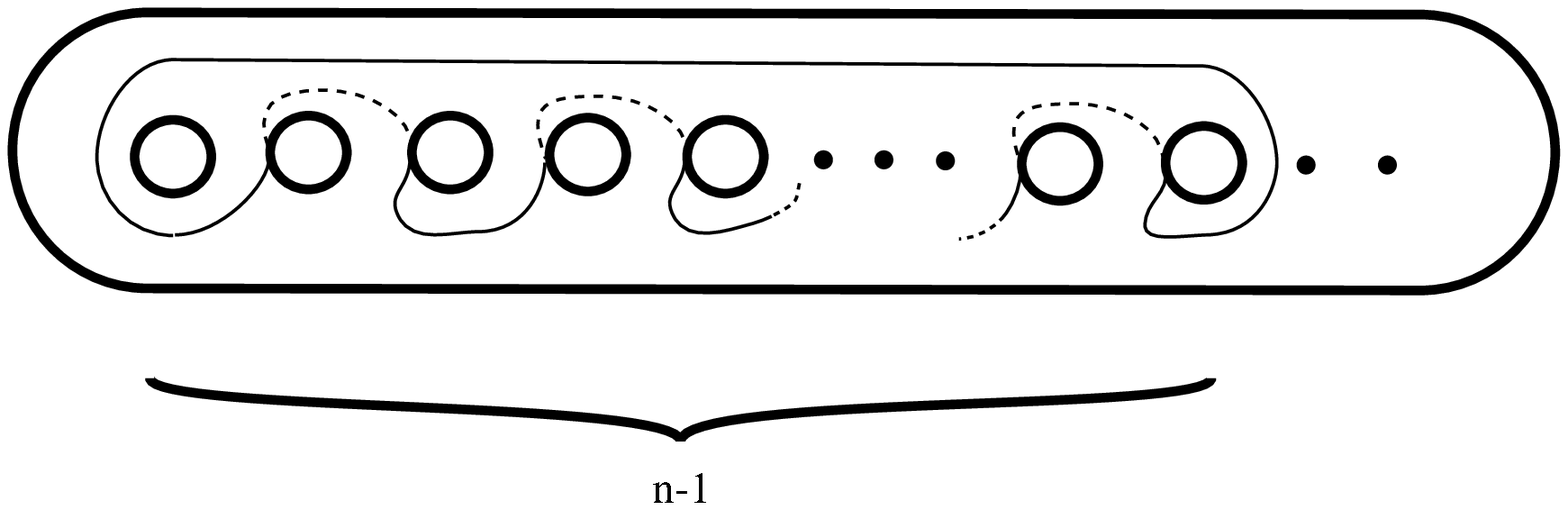}
\nocolon
\caption{}
\label{fig:wave0}
\end{figure}
Figure \ref{fig:wave0} illustrates 
$u = \Prod{l=n-1}{1}{\Prod{i=2l}{n+l-1}{\inv{C_i}}} \cdot
B_n \cdot \Prod{i=n}{2}{C_i} (c_1)$. 
We investigate the action of elements of $G_g$ on $u$. 
\begin{figure}[ht!]
\centering
\includegraphics[height=4cm]{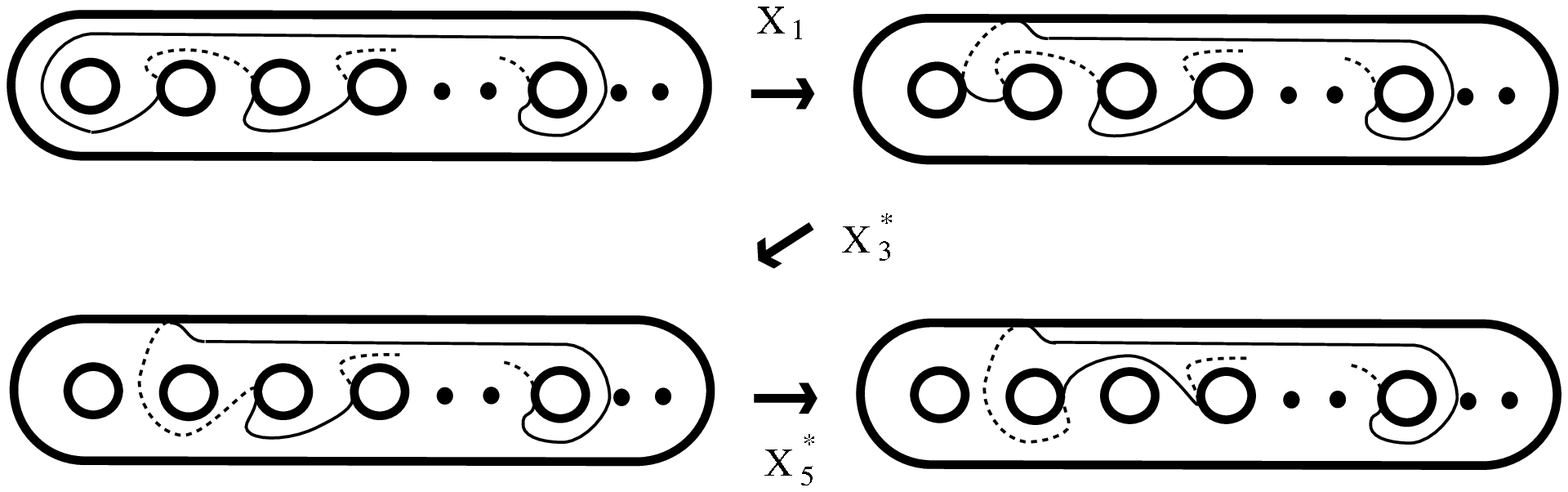}
\nocolon
\caption{}
\label{fig:wave1}
\end{figure}
As indicated in Figure \ref{fig:wave1}, 
$X^*_5 X^*_3 X_1$ acts on $u$. 
We make $\Prod{i=\frac{n}{2} - 2}{2}{X^*_{4i+1} X^*_{4i-1}}$ act 
on this circle.
\begin{figure}[ht!]
\centering
\includegraphics[height=4.5cm]{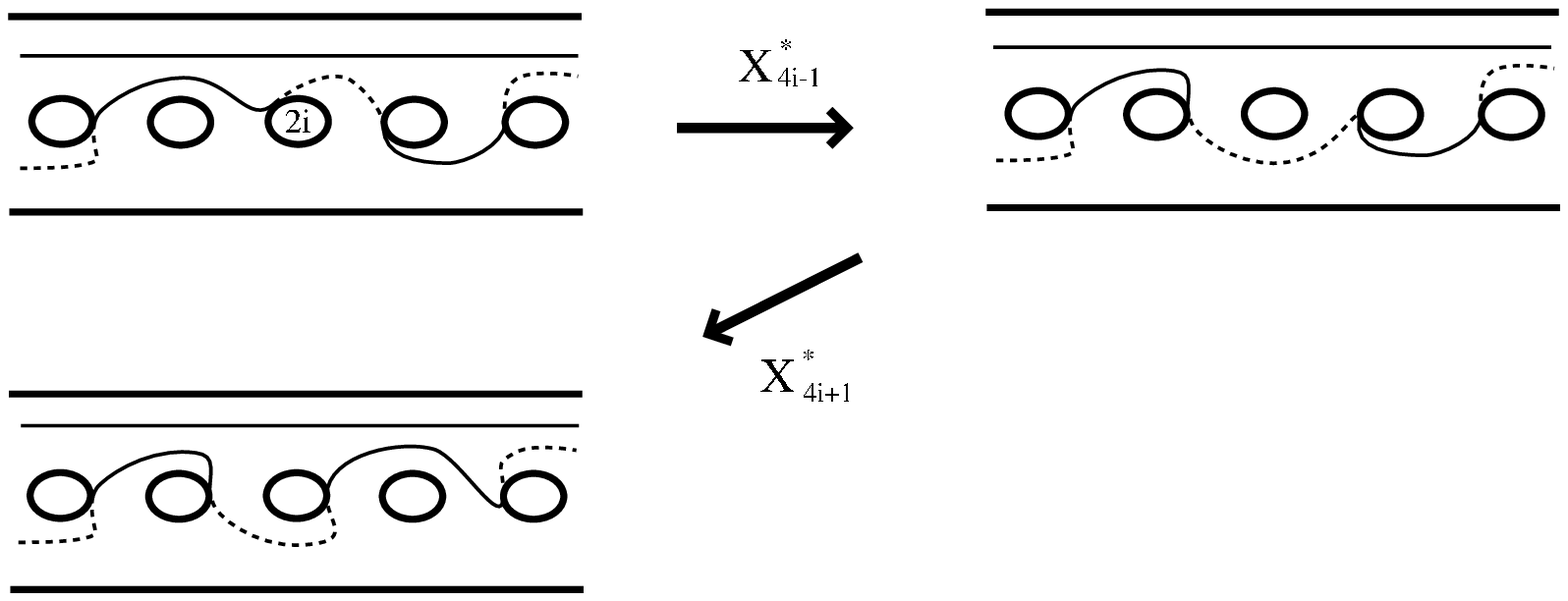}
\nocolon
\caption{}
\label{fig:wave2}
\end{figure} 
In the middle of this action, 
$X^*_{4i+1} X^*_{4i-1}$ acts locally as in Figure \ref{fig:wave2}. 
\begin{figure}[ht!]
\centering
\includegraphics[height=5.5cm]{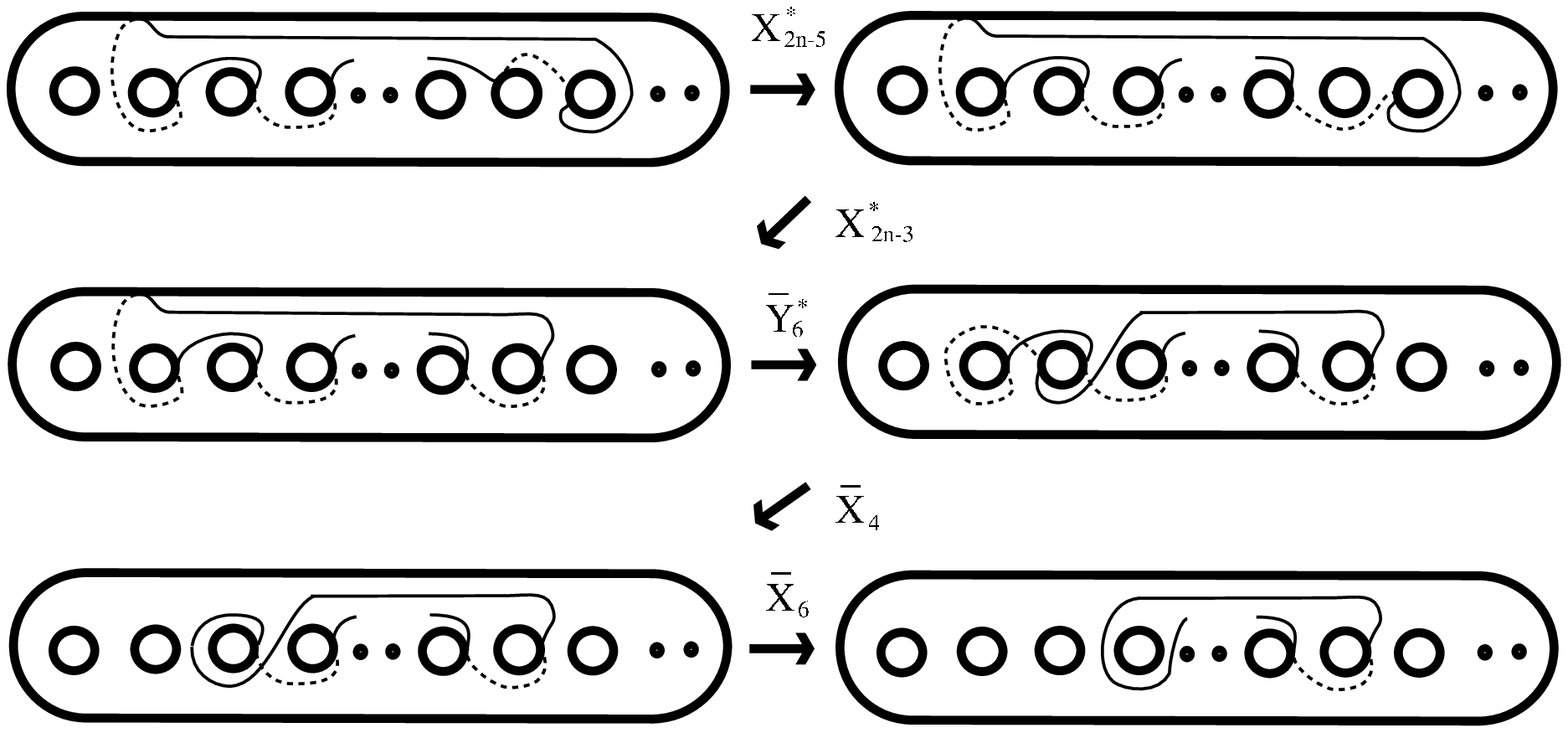}
\nocolon
\caption{}
\label{fig:wave3}
\end{figure} 
Hence, 
$\Prod{i=\frac{n}{2} - 2}{2}{X^*_{4i+1} X^*_{4i-1}} \cdot 
X_1 (u)$ is as the first of Figure \ref{fig:wave3}. 
This figure shows that, by the action of 
$\inv{X_6} \inv{X_4} \inv{Y^*_6} X^*_{2n-3} X^*_{2n-5}$, 
this curve is changed to the $u$ of $n-4$. 
Therefore, for our purpose, it suffices to show that 
$T_u T_u$ is an element of $G_g$ only for $n=4$ or $n=6$. 
\begin{figure}[ht!]
\centering
\includegraphics[height=4.2cm]{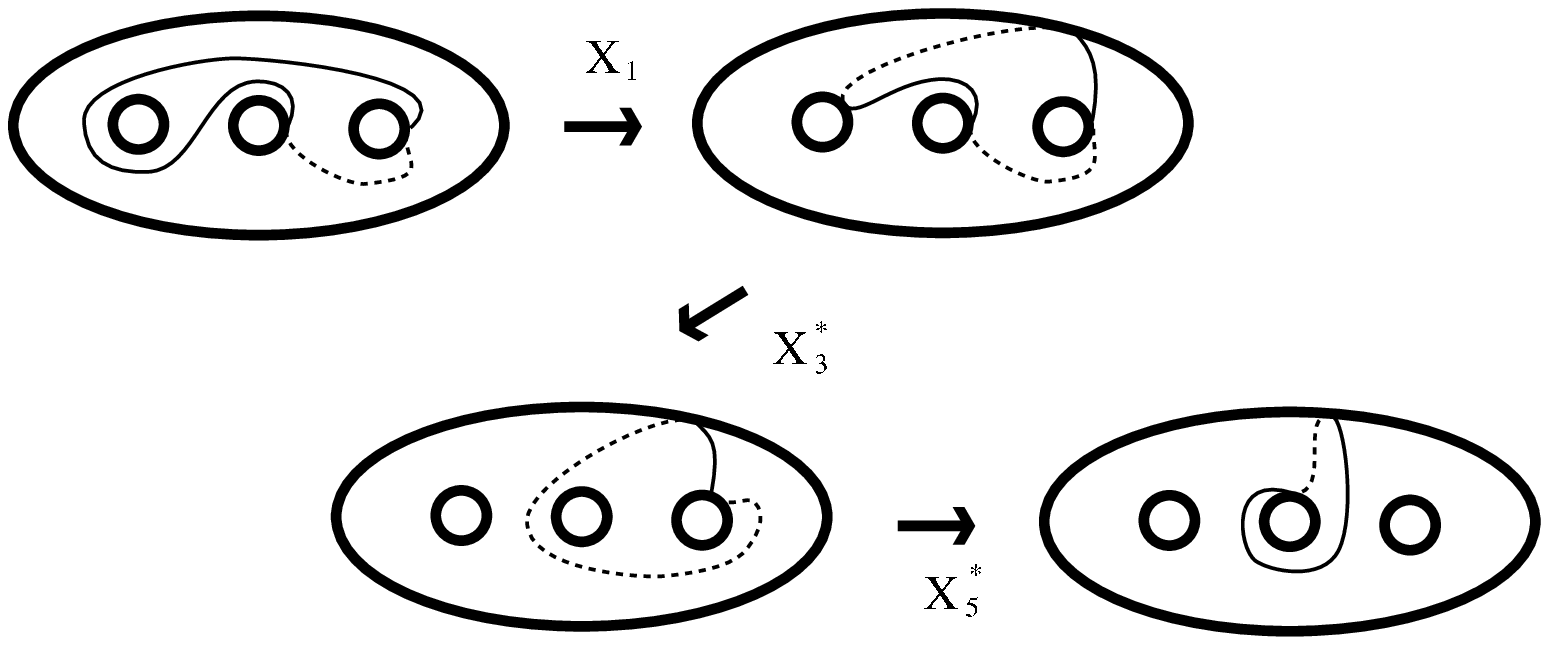}
\nocolon
\caption{}
\label{fig:wave4}
\end{figure} 
Figure \ref{fig:wave4} shows that, when $n=4$, 
$T_u T_u = (\inv{X_1} \inv{X^*_3} \inv{X^*_5})*(Y^*_4 Y^*_4)$. 
\begin{figure}[ht!]
\centering
\includegraphics[height=5cm]{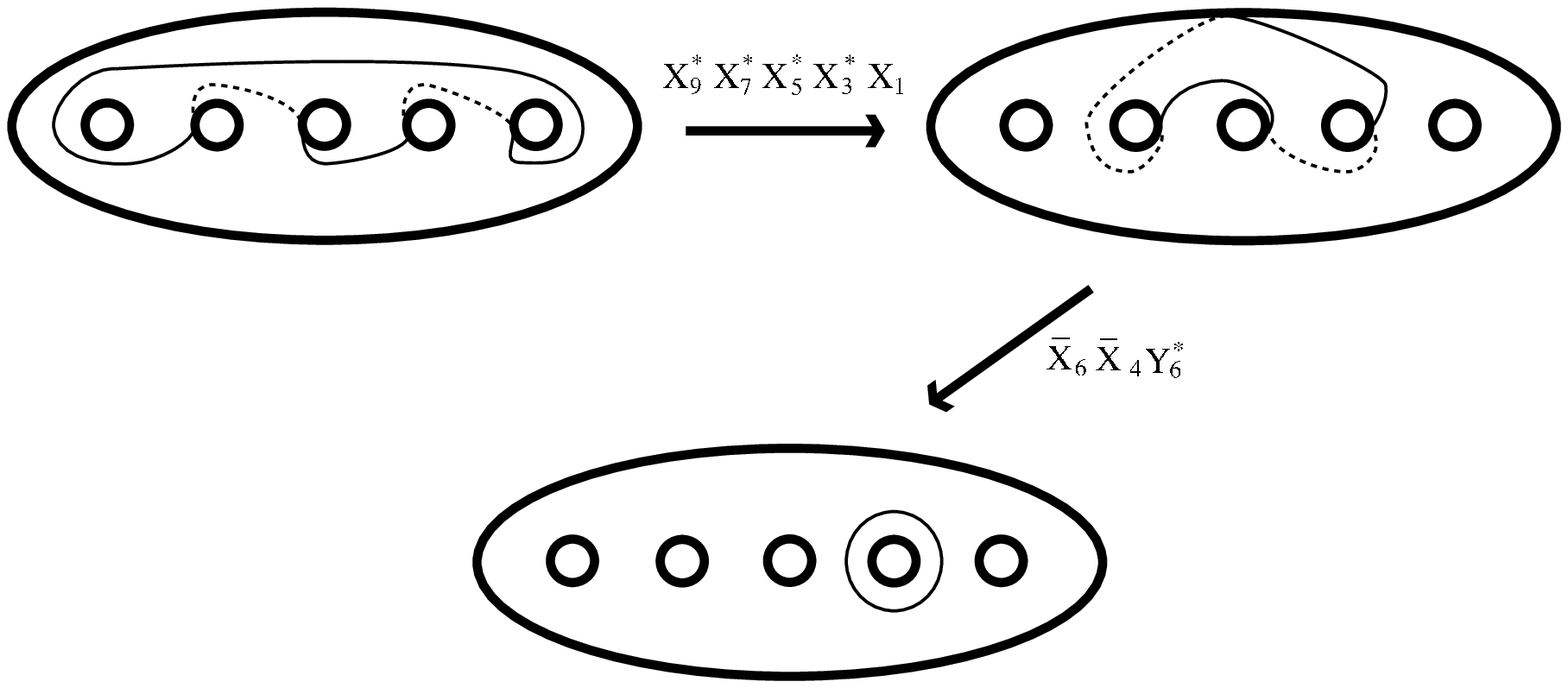}
\nocolon
\caption{}
\label{fig:wave5}
\end{figure}
Figure \ref{fig:wave5} shows that, when $n=5$, 
$T_u T_u = (\inv{X_1} \inv{X^*_3} \inv{X^*_5} \inv{X^*_7} \inv{X^*_9} 
Y^*_6 X_4 X_6) * D_8$. 
\newline
(3) {\it $[1,2,4,6,\ldots,2i,\ldots,2n-2]$ 
$(n \text{ is even, and } 4\leq n \leq g+2)$ are elements of 
$G_g$\/}: 
By (b) of Lemma \ref{lem:Johnson-action}, 
$$
[1,2,4,6,8, \ldots, 2n-2] = 
(\Prod{k=n-2}{1}{\Prod{i=2k+1}{n+k-1}{\inv{C_i}}})*[1,2,3,4,\ldots,n]. 
$$ 
In the same way as (2), 
\begin{equation}
\begin{align*}
[1,2,4,6,8, \ldots, 2n-2] &= 
\Prod{k=2}{n}{ \{ (\Prod{l=n-2}{1}{\Prod{i=2l+1}{n+l-1}{\inv{C_i}}} 
\cdot B_n \cdot \Prod{i=n}{k}{C_i})*(C_{k-1} C_{k-1}) \}}\cdot \\
& \cdot (\Prod{l=n-2}{1}{\Prod{i=2l+1}{n+l-1}{\inv{C_i}}} \cdot 
B_n)*(C_n C_n) \cdot \\
& \cdot \Prod{k=2}{n}{ \{(\Prod{l=n-2}{1}{\Prod{i=2l+1}{n+l-1}{\inv{C_i}}}
\cdot 
\Prod{i=n}{k}{\inv{C_i}} ) * (\inv{C_{k-1}} \inv{C_{k-1}}) \}} \cdot \\ 
& \cdot (\Prod{l=n+2}{1}{\Prod{i=2l+1}{n+l-1}{\inv{C_i}}})*
(\inv{C_n} \inv{C_n}).  
\end{align*}
\end{equation}
By Lemma \ref{lem:Braid-Conj}, 
$(\Prod{l=n-2}{1}{\Prod{i=2l+1}{n+l-1}{\inv{C_i}}}
\cdot 
\Prod{i=n}{k}{\inv{C_i}} ) * (\inv{C_{k-1}} \inv{C_{k-1}})$
and 
\linebreak
$(\Prod{l=n+2}{1}{\Prod{i=2l+1}{n+l-1}{\inv{C_i}}})*
(\inv{C_n} \inv{C_n})$ are elements of $G_g$. 
By the same method as in (2), but using 
$$
\Prod{l=n-2}{1}{\Prod{i=2l+1}{n+l-1}{\inv{C_i}}} \cdot 
\Prod{j=k-2}{1}{C_j} = 
\Prod{j=k-2}{2}{(\inv{C_{2j-1}} C_{2j-2} C_{2j-1})} \cdot 
C_1 \cdot
\Prod{l=n-2}{1}{\Prod{i=2l+1}{n+l-1}{\inv{C_i}}}, 
$$
in place of, 
$$
\Prod{l=n-1}{1}{\Prod{i=2l}{n+l-1}{\inv{C_i}}} \cdot 
\Prod{j=k-2}{1}{C_j} = 
\Prod{j=k-2}{1}{(\inv{C_{2j}} C_{2j-1} C_{2j})} \cdot
\Prod{l=n-1}{1}{\Prod{i=2l}{n+l-1}{\inv{C_i}}},    
$$
we conclude that, for our purpose, 
it suffices to show that 
$ (\Prod{l=n-2}{1}{\Prod{i=2l+1}{n+l-1}{\inv{C_i}}} \cdot B_n \cdot 
\Prod{i=n}{2}{C_i})*(C_i C_i)$ and 
$ (C_1 \Prod{l=n-2}{1}{\Prod{i=2l+1}{n+l-1}{\inv{C_i}}} 
\cdot B_n \cdot \Prod{i=n}{2}{C_i})*(C_i C_i)$ are 
elements of $G_g$. 
\begin{figure}[ht!]
\centering
\includegraphics[height=6cm]{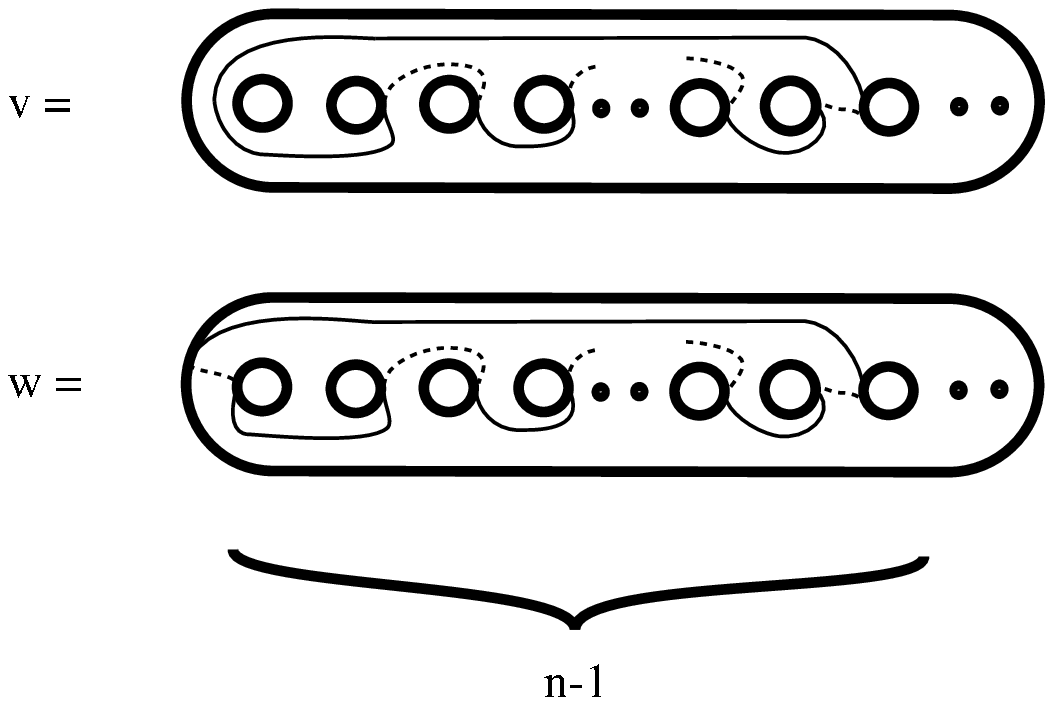}
\nocolon
\caption{}
\label{fig:wavy0}
\end{figure}
Figure \ref{fig:wavy0} illustrates 
$v = \Prod{l=n-2}{1}{\Prod{i=2l+1}{n+l-1}{\inv{C_i}}} \cdot B_n \cdot 
\Prod{i=n}{2}{C_i}(c_1)$ and $w = C_1(v)$. 
First we investigate the actions of elements of $G_g$ on 
$v$. 
In the following argument, we will refer the pictures in 
Figure \ref{fig:wavy1} and Figure \ref{fig:wavy4} by the number 
with (). 
\begin{figure}[htp!]
\centering
\includegraphics[height=18cm]{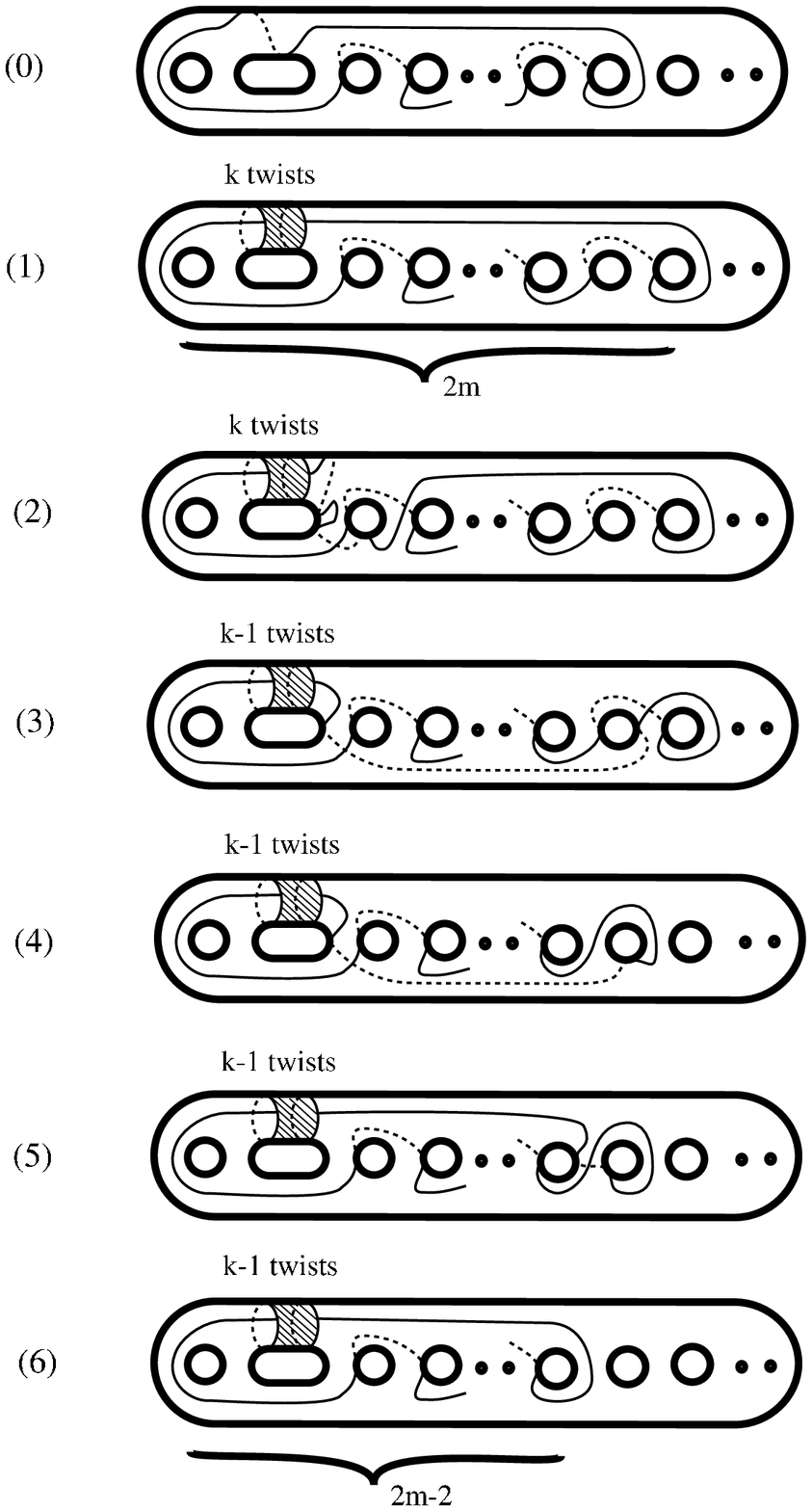}
\nocolon
\caption{}
\label{fig:wavy1}
\end{figure}
By the action of $T_2 \inv{DB_2}$, $v$ is changed to (0). 
Now, we show (1) is $G_g$-equivalent to (6). 
(1) is altered to (2) by the action of $Y^*_6$. 
\begin{figure}[ht!]
\centering
\includegraphics[height=1.3cm]{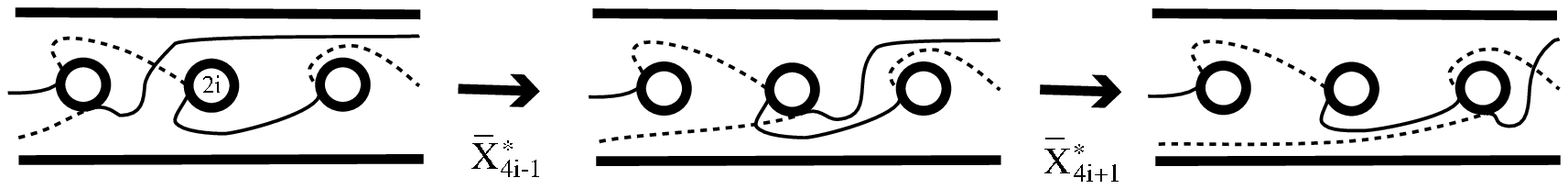}
\nocolon
\caption{}
\label{fig:wavy2}
\end{figure}
We make a sequence of $\inv{X^*_{4i+1}} \inv{X^*_{4i-1}}$'s act 
on this circle. In the middle of this process, 
each $\inv{X^*_{4i+1}} \inv{X^*_{4i-1}}$ acts locally as indicated 
in Figure \ref{fig:wavy2}. 
Hence, (2) is $G_g$-equivalent to (3). 
By the action of $\inv{X^*_{4m-1}}$, (3) is deformed to (4). 
\begin{figure}[ht!]
\centering
\includegraphics[height=1.2cm]{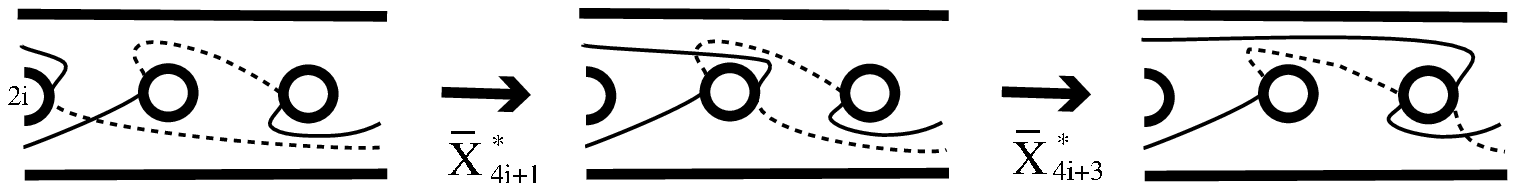}
\nocolon
\caption{}
\label{fig:wavy3}
\end{figure}
In the middle of a sequential action of 
$\inv{X^*_{4i+3}} \inv{X^*_{4i+1}}$'s, 
each $\inv{X^*_{4i+3}} \inv{X^*_{4i+1}}$ acts locally as shown in 
Figure \ref{fig:wavy3}. Hence, 
(4) and (5) are $G_g$-equivalent. 
As a result of the action of $\inv{X^*_{4m-3}}$, 
(5) is altered to (6). 
The above argument shows that (1) is $G_g$-equivalent to (6). 
\begin{figure}[ht!]
\centering
\includegraphics[height=6.4cm]{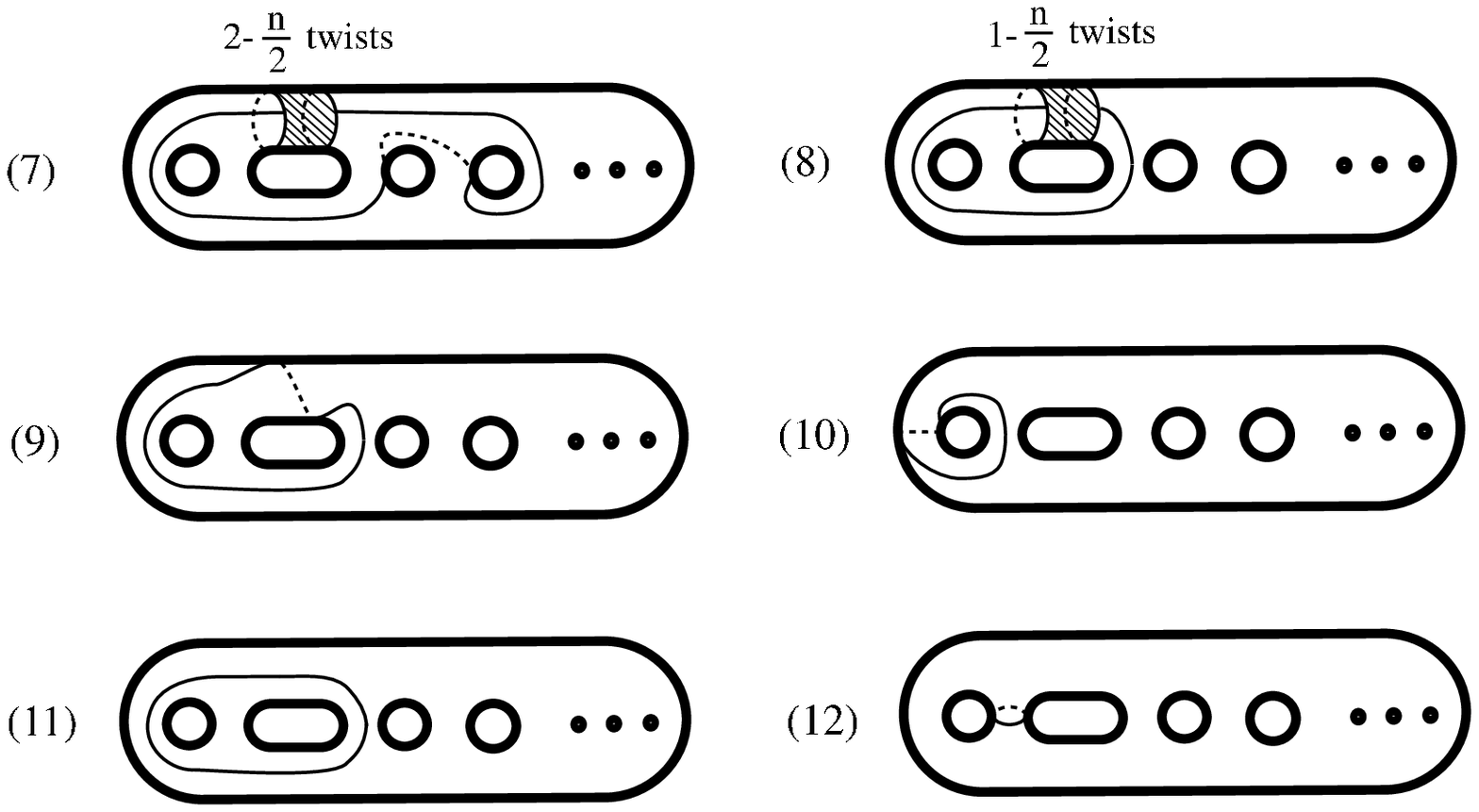}
\nocolon
\caption{}
\label{fig:wavy4}
\end{figure}
For (0), we apply the above process from (1) to (6) repeatedly, 
then we get (7). 
The element $\inv{X^*_5} \inv{X^*_7} \inv{Y^*_6}$ alters (7) into (8). 
If $\frac{n}{2}$ is even, $DB_4^{\frac{n}{4}-1}$ deforms (8) into (9). 
Since (9) is changed to (10) by the action of $\inv{X_3}$, 
there exists an element $h$ of $G_g$ such that 
$h*(T_v T_v) = X_1 X_1$. 
If $\frac{n}{2}$ is odd, $DB_4^{\frac{n-2}{4}}$ deforms (8) into (11). 
Since (11) is changed to (12) by the action of $X_1 \inv{Y^*_4}$, 
there exists an element $h$ of $G_g$ such that 
$h*(T_v T_v) = D_3$. 
Next, we investigate the actions of $G_g$ on $w$. 
\begin{figure}[ht!]
\centering
\includegraphics[height=4.2cm]{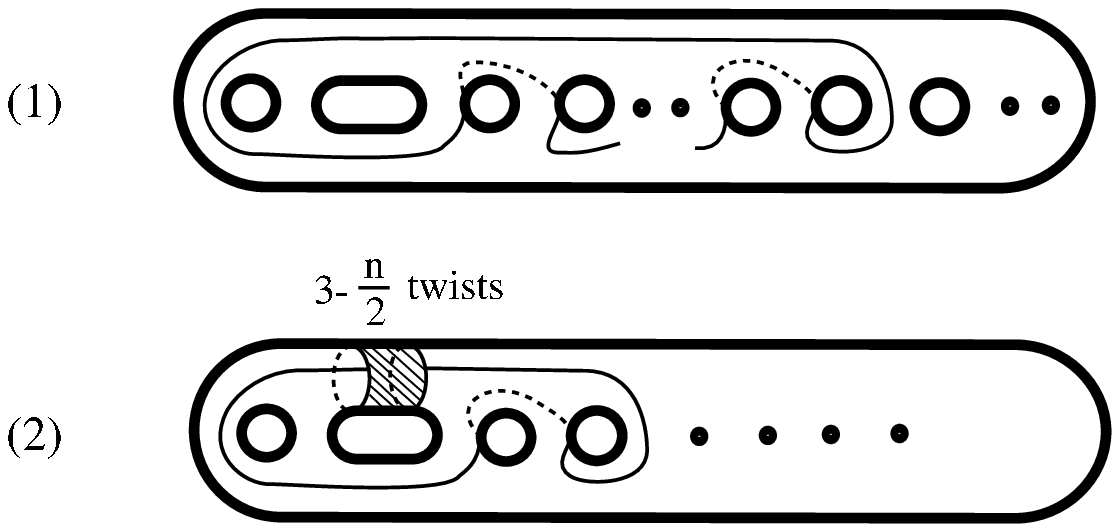}
\nocolon
\caption{}
\label{fig:wavy5}
\end{figure}
The action of $\inv{T_1} T_2$ deforms $w$ into (1) 
of Figure \ref{fig:wavy5}. 
After the repeated application of the actions from (1) to (6) of Figure
\ref{fig:wavy1},  this circle is altered to (2) of Figure \ref{fig:wavy5}. 
By the same argument for $v$, 
when $\frac{n}{2}$ is even, there is a $h$ of $G_g$ such that 
$h*(T_w T_w) = D_3$, on the other hand, 
when $\frac{n}{2}$ is odd, there is a $h$ of $G_g$ such that 
$h*(T_w T_w) = X_1 X_1$. 
Therefore, $[1,2,4,6,8,\ldots,2n-2]$ is an element of 
$G_g$.  
\end{proof}
We prove that any odd subchain map of 
$(c_1, c_2, c_3, \ldots, c_{2g+1})$ or 
$(c_{\beta}, c_5, c_6, \ldots,
$ $
 c_{2g})$ is a product of 
elements listed on Lemma \ref{lem:special-elements} and 
elements of $G_g$. 
The following lemma shows that 
any odd subchain map of $(c_{\beta}, c_5, c_6, \ldots, c_{2g})$ 
is a product of an odd subchain map of 
$(c_1, c_2, c_3, \ldots, c_{2g+1})$ and elements of $G_g$. 
\begin{lem}\label{lem:beta-to-straight}
%
$D_3 \inv{T_1} (c_{\beta}) = c_3 + c_4$. 
\end{lem}
\begin{proof}
Figure \ref{fig:BtoS} proves this lemma.
\end{proof}
From here to the end of this subsection, 
odd subchain maps mean only those of $(c_1, c_2, c_3, \ldots, c_{2g+1})$. 
The following lemma shows that any odd subchain map, whose length is 
at least 5 and which begins from $1,2,3,4$, is a product 
of shorter odd subchain maps and elements of $G_g$. 
\begin{lem}\label{lem:chain-shortening}
%
\begin{equation}
\begin{align*}
[1,2,3,4] [1,2,3,5]^{-1} &[1,2,3,4][1,2,4,6,7,\ldots,2n] \cdot\\
\cdot
(C_4 B_4 \inv{C_4})*[3,4,5,\ldots,2n] 
&= [4,6,7,\ldots,2n][1,2,3,4,\ldots,2n]
\end{align*}
\end{equation}
\end{lem}
\begin{proof}
By (a) of Lemma \ref{lem:Johnson-action}, 
$\inv{C_4}*[3,4,5,\ldots,2n] = [3,4,5,\ldots,2n]$, and by (d) of Lemma 
\ref{lem:Johnson-action}, 
\begin{align*}
[1,2,3,4][1,2,5,6, \ldots, 2n] &\cdot (B_4 \inv{C_4})*[3,4,5,\ldots,2n] = \\
&= [5,6, \ldots, 2n][1,2,3,4,\ldots,2n].
\end{align*}
By applying $C_4$ to the above equation, we get the equation which we need. 
\end{proof}
For any odd subchain map $[i_1, i_2, \ldots, i_r]$, we define a sequence 
$[[\epsilon_1, \epsilon_2, \ldots, 
$ $
\epsilon_{2g+2}]]$ as follows:
$\epsilon_k =1$ if $k$ is a member of $\{ i_1, i_2, \ldots, i_r\}$, 
and $\epsilon_k = 0$ if $k$ is not a member of 
$\{ i_1, i_2, \ldots, i_r\}$. 
For this sequence $[[\epsilon_1, \epsilon_2, \ldots, \epsilon_{2g+2}]]$, 
we construct the sequence $[[\delta_1, \delta_2, \ldots, \delta_{2g+2}]]$ 
by the following rule: 
$(\delta_{2i-1}, \delta_{2i}) = (0,0)$ if 
$(\epsilon_{2i-1}, \epsilon_{2i})=(0,0)$, 
$(\delta_{2i-1}, \delta_{2i}) = (1,0)$ if 
$(\epsilon_{2i-1}, \epsilon_{2i})=(0,1)$,  
$(\delta_{2i-1}, \delta_{2i}) = (0,1)$ if 
$(\epsilon_{2i-1}, \epsilon_{2i})=(1,0)$, 
$(\delta_{2i-1}, \delta_{2i}) = (1,1)$ if 
$(\epsilon_{2i-1}, \epsilon_{2i})=(1,1)$. 
The odd subchain map $[j_1,j_2,\ldots,j_r]$, which corresponds to 
the sequence $[[\delta_1, \delta_2, \ldots, 
$ $
\delta_{2g+2}]]$, 
is called {\it the reversion\/} of 
$[i_1, i_2, \ldots, i_r]$. 
\begin{lem}\label{lem:Gg-action}
%
{\rm(1)}\qua For any odd subchain map $c$, there is an element of $G_g$ which brings 
$c$ to its reversion. \newline
{\rm(2)}\qua When $k \leq i-3$, 
$(\inv{C_{i-1}} C_{i-2} C_{i-1})*[\ldots,k,i,j,\ldots] 
=[\ldots,k,i-2,j,\ldots]$. \newline
{\rm(3)}\qua When $k \leq i-2$, 
$(C_i C_{i-1} \inv{C_i})*[\ldots,k,i,i+1,\ldots]
=[\ldots,k,i-1,i,\ldots]$. 
\end{lem}
\begin{proof}
Lemma \ref{lem:Johnson-action} shows (2) and (3). 
Since, $\inv{T_1} T_2 = \inv{C_1} \inv{C_3} C_5 \cdots C_{2g+1}$ 
and $D_{2i-1} = C_{2i-1} C_{2i-1}$ $(1 \leq i \leq g+1)$ are 
elements of $G_g$, 
$C_1^{\pm 1} C_3^{\pm 1} \cdots C_{2g+1}^{\pm 1}$ 
is an elements of $G_g$ for any choice of $\pm 1$'s. 
Let $[[\epsilon_1, \epsilon_2, \ldots, \epsilon_{2g+2}]]$ 
be the 0-1 sequence corresponding to 
$[i_1, i_2, \ldots, i_r]$. 
We define $\gamma_i$ $(1 \leq i \leq g+1)$ as follows: 
$\gamma_i = +1$ if 
$(\epsilon_{2i-1}, \epsilon_{2i}) = (0,0), (0,1)$, or $(1,1)$, 
and 
$\gamma_i = -1$ if 
$(\epsilon_{2i-1}, \epsilon_{2i}) = (1,0)$. 
Then 
$(C_1^{\gamma_1} C_3^{\gamma_2} \cdots C_{2g+1}^{\gamma_{g+1}})
*[i_1,\ldots,i_r]$ is the reversion of $[i_1, \ldots, i_r]$. 
\end{proof}   
By (2) of the above lemma, any odd subchain map is 
deformed to an odd subchain map $[i_1, i_2, \ldots, i_r]$ 
such that $i_{l+1} - i_l \leq 2$ under the action of $G_g$. 
If there are at least two disjoint pairs of indices $(i_l, i_{l+1})$ in an
odd subchain map $[i_1, i_2, \ldots, i_r]$ such that 
$i_{l+1} = i_l + 1$, then, 
by (3) of the above lemma, 
this odd subchain map is altered to the odd subchain map 
which begins from $1,2,3,4$ under the action of $G_g$. 
Therefore, by Lemma \ref{lem:chain-shortening}, this odd subchain map 
is a product of shorter odd subchain maps and elements of $G_g$. 
Hence, it suffices to show that $[1,3,5,7,9, \ldots]$, 
$[2,4,6,8,10, \ldots]$, $[1,2,3,5,7, \ldots]$, $[1,2,4,6,8,\ldots]$, 
and $[1,2,3,4]$ are elements of $G_g$. 
By (1) of Lemma \ref{lem:Gg-action}, the second ones are changed to 
the first ones, and the third ones are changed to 
the fourth ones by the action of $G_g$. 
On the other hand, we have already shown that 
$[1,3,5,7,9,\ldots]$, $[1,2,4,6,8,\ldots]$, and $[1,2,3,4]$ are 
elements of $G_g$ in Lemma \ref{lem:special-elements}. 
Therefore, Lemma \ref{lem:Torelli} is proved. 
\subsection{The level 2 prime congruence subgroup of 
Sp $(2g, {\Bbb Z})$}
In this subsection, we assume $g \geq 3$. 
Let $\Phi_2$ be the natural homomorphism from ${\cal M}_g$ to 
$\Symp (2g, {\Bbb Z}_2)$ defined by the action of ${\cal M}_g$ 
on the ${\Bbb Z}_2$-coefficient first homology group 
$H_1(\Sigma_g ; {\Bbb Z}_2)$. 
In this section, we show the following lemma. 
\begin{lem}\label{lem:level2-congruence}
%
$\ker \Phi_2$ is a subgroup of $G_g$. 
\end{lem}
We denote the kernel of the natural homomorphism from 
$\Symp(2g,{\Bbb Z})$ to $\Symp(2g,{\Bbb Z}_2)$ by $\Symp^{(2)} (2g)$. 
We set a basis of $H_1(\Sigma_g; {\Bbb Z})$ as in Figure \ref{fig:basis}, 
and define the intersection form $(,)$ on $H_1(\Sigma_g; {\Bbb Z})$ 
to satisfy $(x_i, y_j)=\delta_{i,j}$, $(x_i, x_j) = (y_i, y_j) = 0$ 
$(1 \leq i, j, \leq g)$. 
An element $a$ of $H_1(\Sigma_g; {\Bbb Z})$ is called 
{\it primitive\/} if there is no element $n (\not=0, \pm 1)$ of 
${\Bbb Z}$, and no element $b$ of $H_1(\Sigma_g; {\Bbb Z})$ such that 
$a = nb$. 
For a primitive element $a$ of $H_1(\Sigma_g; {\Bbb Z})$, we define 
an isomorphism $T_a\co H_1(\Sigma_g; {\Bbb Z}) \to H_1(\Sigma_g; {\Bbb Z})$ 
by $T_a(v) = v + (a,v)a$. 
This isomorphism is the same as the action of Dehn twist about a simple
closed curve  representing $a$ on $H_1(\Sigma_g; {\Bbb Z})$.  
We call $T_a^2$ {\it the square transvection about $a$\/}. 
Johnson \cite{Johnson2} showed the following result. 
\begin{lem}\label{lem:double-transvection}
%
$\Symp^{(2)} (2g)$ is generated by square transvections. 
\end{lem}
$\Symp^{(2)} (2g)$ is finitely generated. In fact, we show: 
\begin{lem}\label{lem:level2-congruence-generator}
%
$\Symp^{(2)} (2g)$ is generated by the square transvections about 
the primitive elements 
$\sum_{i=1}^g (\epsilon_i x_i + \delta_i y_i)$, 
where $\epsilon_i =0,1$ and $\delta_i =0,1$. 
\end{lem}
We define, for any primitive element $a$ and $b$ of 
$H_1(\Sigma_g; {\Bbb Z})$, two operation $\boxplus$ and 
$\boxminus$ by 
$$
a \boxplus b = a+2(a,b)b,\quad  a \boxminus b = a-2(a,b)b.
$$ 
We remark that $T_{a\boxplus b}^2 = T_b^{-2} T_a^2 T_b^2$, 
$T_{a\boxminus b}^2 = T_b^2 T_a^2 T_b^{-2}$, and 
$(a \boxplus b)\boxminus b = a = (a \boxminus b) \boxplus b$. 
We denote the element $\sum_{i=1}^g (a_i^1 x_i + a_i^2 y_i)$ 
of $H_1(\Sigma_g; {\Bbb Z})$ by 
$[(a_1^1, a_1^2),(a_2^1, a_2^2), \cdots, (a_g^1, a_g^2)]$, 
and call each $(a_i^1, a_i^2)$ as {\it a block\/}. 
For a positive integer $k$, $a (\boxplus b)^k$ is the result of
the $k$-fold application of $\boxplus b$ on $a$, and 
$a (\boxplus b)^{-k}$ is the result of the $k$-fold 
application of $\boxminus b$ on $a$. 
\begin{lem}\label{lem:first-Euclid}
%
For any primitive element $a$ of $H_1(\Sigma_g; {\Bbb Z})$, 
by applying 
$\boxplus [(0,0), 
$ $
\ldots, (0,0), (1,0), (0,0), \ldots, (0,0)]$ or 
$\boxplus [(0,0), \ldots, (0,0), (0,1), (0,0), \ldots, (0,0)]$ 
several times, 
each block of $a$ is altered to $(0,0)$, $(p,0)$, $(0,p)$, or $(p,p)$. 
\end{lem}
\begin{proof}
Let $(m,n)$ be the $i$-th block of $a$. 
First we consider the case when $|m| > |n| \not= 0 $.
There is an integer $k$ such that $|m-2kn| \leq |n|$. 
Let $e_i$ be the element of $H_1(\Sigma_g; {\Bbb Z})$, 
the $i$-th block of which is $(1,0)$, and every other block of which
is $(0,0)$. 
Since, 
$ [\cdots, (m,n), \cdots] \boxplus e_i = [\cdots, (m-2n, n), \cdots], $ 
and 
$[\cdots, (m,n), \cdots] \boxminus e_i = [\cdots, (m+2n,n), \cdots],$
we get  
$[\cdots, (m,n), \cdots] (\boxplus e_i)^k = 
[\cdots, (m-2kn,n),\cdots]$. 
This means that, by repeated application of $\boxplus e_i$, 
the $i$-th block $(m,n)$ is altered such that 
$|m| \leq |n|$. 
Next, we consider the case when $0 \not= |m| < |n|$.   
Let $f_i$ be the element of $H_1(\Sigma_g; {\Bbb Z})$, 
the $i$-th block of which is $(0,1)$, and other blocks of which 
are $(0,0)$. 
Since, 
$ [\cdots, (m,n), \cdots] \boxplus f_i = [\cdots, (m, n+2m), \cdots], $ 
and 
$[\cdots, (m,n), \cdots] \boxminus f_i = [\cdots, (m,n-2m), \cdots],$
by the same argument as the previous case, 
by repeated application of $\boxplus f_i$, 
the $i$-th block is altered such that $|m| \geq |n|$. 
The above arguments show that, 
after several application of $\boxplus e_i$ or $\boxplus f_i$, 
the $i$-th block $(m,n)$ of $a$ is altered to be  
$|m|=|n|$, or $m=0$, or $n=0$. 
If $n=-m$, the $i$-th block changed to $(m,m)$ by the 
application of $\boxplus f_i$. 
For each $i$-th block, we do the same operation as above. 
Then, this lemma follows.   
\end{proof}
For a primitive element of $H_1(\Sigma_g; {\Bbb Z})$, each of whose blocks
is $(p,0)$, or $(0,p)$, or $(p,p)$, (where $p$ can be different from block
to block) we apply several operations 
$\boxplus [\ldots, (\epsilon_i, \delta_i), \ldots]$, 
where $\epsilon_i = 0,1$ and $\delta_i =0,1$. 
Then we obtain the following equations, where $\Ccirc$ means a sequence 
of $(0,0)$, and $\Cbullet$ means the part which is not changed. 
{\allowdisplaybreaks
\begin{align*}
[\Cbullet, &(p,0), (q,0), \Cbullet] \boxminus [\Ccirc, (1,0),(0,1), \Ccirc] 
\boxplus [\Ccirc, (0,0), (0,1), \Ccirc] \\
&=[\Cbullet, (p-2q,0),(q,0), \Cbullet], \\
[\Cbullet, &(p,0), (q,0), \Cbullet] \boxplus [\Ccirc, (1,0),(0,1), \Ccirc] 
\boxminus [\Ccirc, (0,0), (0,1), \Ccirc] \\
&=[\Cbullet, (p+2q,0),(q,0), \Cbullet], \\
[\Cbullet, &(p,0), (q,0), \Cbullet] \boxminus [\Ccirc, (0,1),(1,0), \Ccirc] 
\boxplus [\Ccirc, (0,1), (0,0), \Ccirc] \\
&=[\Cbullet, (p,0),(q-2p,0), \Cbullet], \\
[\Cbullet, &(p,0), (q,0), \Cbullet] \boxplus [\Ccirc,(0,1),(1,0),\Ccirc] 
\boxplus [\Ccirc, (0,1), (0,0), \Ccirc] \\
&=[\Cbullet, (p,0),(q+2p,0), \Cbullet], 
\end{align*}
}
%
%
{\allowdisplaybreaks
\begin{align*}
[\Cbullet, &(0,p), (0,q), \Cbullet] \boxplus [\Ccirc, (0,1),(1,0),\Ccirc] 
\boxminus [\Ccirc, (0,0), (1,0), \Ccirc] \\
&=[\Cbullet, (0,p-2q),(0,q), \Cbullet], \\
[\Cbullet, &(0,p), (0,q), \Cbullet] \boxminus [\Ccirc, (0,1),(1,0), \Ccirc] 
\boxplus [\Ccirc, (0,0), (1,0), \Ccirc] \\
&=[\Cbullet, (0,p+2q),(0,q), \Cbullet], \\
[\Cbullet, &(0,p), (0,q), \Cbullet] \boxplus [\Ccirc, (1,0),(0,1),\Ccirc] 
\boxminus [\Ccirc, (1,0), (0,0), \Ccirc] \\
&=[\Cbullet, (0,p),(0,q-2p), \Cbullet], \\
[\Cbullet, &(0,p), (0,q), \Cbullet] \boxminus [\Ccirc,(1,0),(0,1),\Ccirc] 
\boxplus [\Ccirc, (1,0), (0,0), \Ccirc] \\
&=[\Cbullet, (0,p),(0,q+2p), \Cbullet], 
\end{align*}
}
%
%
{\allowdisplaybreaks
\begin{align*}
[\Cbullet, &(p,0), (0,q), \Cbullet] \boxplus [\Ccirc, (1,0),(1,0),\Ccirc] 
\boxminus [\Ccirc, (0,0), (1,0), \Ccirc] \\
&=[\Cbullet, (p-2q,0),(0,q), \Cbullet], \\
[\Cbullet, &(p,0), (0,q), \Cbullet] \boxminus [\Ccirc, (1,0),(1,0),\Ccirc] 
\boxplus [\Ccirc, (0,0), (1,0), \Ccirc] \\
&=[\Cbullet, (p+2q,0),(0,q), \Cbullet], \\
[\Cbullet, &(p,0), (0,q), \Cbullet] \boxminus [\Ccirc,(0,1),(0,1),\Ccirc] 
\boxplus [\Ccirc, (0,1), (0,0), \Ccirc] \\
&=[\Cbullet, (p,0),(0,q-2p), \Cbullet], \\
[\Cbullet, &(p,0), (0,q), \Cbullet] \boxplus [\Ccirc,(0,1),(0,1),\Ccirc] 
\boxplus [\Ccirc, (0,1), (0,0), \Ccirc] \\
&=[\Cbullet, (p,0),(0,q+2p), \Cbullet], 
\end{align*}
}
%
%
{\allowdisplaybreaks
\begin{align*}
[\Cbullet, &(0,p), (q,0), \Cbullet] \boxminus [\Ccirc,(0,1),(0,1),\Ccirc] 
\boxplus [\Ccirc, (0,0), (0,1), \Ccirc] \\
&=[\Cbullet, (0,p-2q),(q,0), \Cbullet], \\
[\Cbullet, &(0,p), (q,0), \Cbullet] \boxplus [\Ccirc,(0,1),(0,1),\Ccirc] 
\boxminus [\Ccirc, (0,0), (0,1), \Ccirc] \\
&=[\Cbullet, (0,p+2q),(q,0), \Cbullet], \\
[\Cbullet, &(0,p), (q,0), \Cbullet] \boxplus [\Ccirc,(1,0),(1,0),\Ccirc] 
\boxminus [\Ccirc, (1,0), (0,0), \Ccirc] \\
&=[\Cbullet, (0,p),(q-2p,0), \Cbullet], \\
[\Cbullet, &(0,p), (q,0), \Cbullet] \boxminus [\Ccirc,(1,0),(1,0),\Ccirc] 
\boxplus [\Ccirc, (1,0), (0,0), \Ccirc] \\
&=[\Cbullet, (0,p),(q+2p,0), \Cbullet], 
\end{align*}
}
%
%
{\allowdisplaybreaks
\begin{align*}
[\Cbullet, &(0,p), (q,q), \Cbullet] \boxminus [\Ccirc,(0,1),(0,1),\Ccirc] 
\boxplus [\Ccirc, (0,0), (0,1), \Ccirc] \\
&=[\Cbullet, (0,p-2q),(q,q), \Cbullet], \\
[\Cbullet, &(0,p), (q,q), \Cbullet] \boxplus [\Ccirc,(0,1),(0,1),\Ccirc] 
\boxminus [\Ccirc, (0,0), (0,1), \Ccirc] \\
&=[\Cbullet, (0,p+2q),(q,q), \Cbullet], \\
[\Cbullet, &(0,p), (q,q), \Cbullet] \boxplus [\Ccirc,(1,0),(1,1),\Ccirc] 
\boxminus [\Ccirc, (1,0), (0,0), \Ccirc] \\
&=[\Cbullet, (0,p),(q-2p, q-2p), \Cbullet], \\
[\Cbullet, &(0,p), (q,q), \Cbullet] \boxminus [\Ccirc,(1,0),(1,1),\Ccirc] 
\boxplus [\Ccirc, (1,0), (0,0), \Ccirc] \\
&=[\Cbullet, (0,p),(q+2p,q+2p), \Cbullet], 
\end{align*}
}
%
%
{\allowdisplaybreaks
\begin{align*}
[\Cbullet, &(p,p), (0,q), \Cbullet] \boxplus [\Ccirc,(1,1),(1,0),\Ccirc] 
\boxminus [\Ccirc, (0,0), (1,0), \Ccirc] \\
&=[\Cbullet, (p-2q,p-2q),(0,q), \Cbullet], \\
[\Cbullet, &(p,p), (0,q), \Cbullet] \boxminus [\Ccirc,(1,1),(1,0),\Ccirc] 
\boxplus [\Ccirc, (0,0), (1,0), \Ccirc] \\
&=[\Cbullet, (p+2q,p+2q),(0,q), \Cbullet], \\
[\Cbullet, &(p,p), (0,q), \Cbullet] \boxminus [\Ccirc,(0,1),(0,1),\Ccirc] 
\boxplus [\Ccirc, (0,1), (0,0), \Ccirc] \\
&=[\Cbullet, (p,p),(0,q-2p), \Cbullet], \\
[\Cbullet, &(p,p), (0,q), \Cbullet] \boxplus [\Ccirc,(0,1),(0,1),\Ccirc] 
\boxminus [\Ccirc, (0,1), (0,0), \Ccirc] \\
&=[\Cbullet, (p,p),(0,q+2p), \Cbullet], 
\end{align*}
}
%
%
{\allowdisplaybreaks
\begin{align*}
[\Cbullet, &(p,0), (q,q), \Cbullet] \boxminus [\Ccirc,(1,0),(0,1),\Ccirc] 
\boxplus [\Ccirc, (0,0), (0,1), \Ccirc] \\
&=[\Cbullet, (p-2q,0),(q,q), \Cbullet], \\
[\Cbullet, &(p,0), (q,q), \Cbullet] \boxplus [\Ccirc,(1,0),(0,1),\Ccirc] 
\boxminus [\Ccirc, (0,0), (0,1), \Ccirc] \\
&=[\Cbullet, (p+2q,0),(q,q), \Cbullet], \\
[\Cbullet, &(p,0), (q,q), \Cbullet] \boxminus [\Ccirc,(0,1),(1,1),\Ccirc] 
\boxplus [\Ccirc, (0,1), (0,0), \Ccirc] \\
&=[\Cbullet, (p,0),(q-2p,q-2p), \Cbullet], \\
[\Cbullet, &(p,0), (q,q), \Cbullet] \boxplus [\Ccirc,(0,1),(1,1),\Ccirc] 
\boxminus [\Ccirc, (0,1), (0,0), \Ccirc] \\
&=[\Cbullet, (p,0),(q+2p,q+2p), \Cbullet], 
\end{align*}
}
%
%
{\allowdisplaybreaks
\begin{align*}
[\Cbullet, &(p,p), (q,0), \Cbullet] \boxminus [\Ccirc,(1,1),(0,1),\Ccirc] 
\boxplus [\Ccirc, (0,0), (0,1), \Ccirc] \\
&=[\Cbullet, (p-2q,p-2q),(q,0), \Cbullet], \\
[\Cbullet, &(p,p), (q,0), \Cbullet] \boxplus [\Ccirc,(1,1),(0,1),\Ccirc] 
\boxminus [\Ccirc, (0,0), (0,1), \Ccirc] \\
&=[\Cbullet, (p+2q,p+2q),(q,0), \Cbullet], \\
[\Cbullet, &(p,p), (q,0), \Cbullet] \boxminus [\Ccirc,(0,1),(1,0),\Ccirc] 
\boxplus [\Ccirc, (0,1), (0,0), \Ccirc] \\
&=[\Cbullet, (p,p),(q-2p,0), \Cbullet], \\
[\Cbullet, &(p,p), (q,0), \Cbullet] \boxplus [\Ccirc,(0,1),(1,0),\Ccirc] 
\boxminus [\Ccirc, (0,1), (0,0), \Ccirc] \\
&=[\Cbullet, (p,p),(q+2p,0), \Cbullet], 
\end{align*}
}
%
%
{\allowdisplaybreaks
\begin{align*}
[\Cbullet, &(p,p), (q,q), \Cbullet] \boxminus [\Ccirc,(1,1),(0,1),\Ccirc] 
\boxplus [\Ccirc, (0,0), (0,1), \Ccirc] \\
&=[\Cbullet, (p-2q,p-2q),(q,q), \Cbullet], \\
[\Cbullet, &(p,p), (q,q), \Cbullet] \boxplus [\Ccirc,(1,1),(0,1),\Ccirc] 
\boxminus [\Ccirc, (0,0), (0,1), \Ccirc] \\
&=[\Cbullet, (p+2q,p+2),(q,q), \Cbullet], \\
[\Cbullet, &(p,p), (q,q), \Cbullet] \boxminus [\Ccirc,(0,1),(1,1),\Ccirc] 
\boxplus [\Ccirc, (0,1), (0,0), \Ccirc] \\
&=[\Cbullet, (p,p),(q-2p,q-2p), \Cbullet], \\
[\Cbullet, &(p,p), (q,q), \Cbullet] \boxplus [\Ccirc,(0,1),(1,1),\Ccirc] 
\boxminus [\Ccirc, (0,1), (0,0), \Ccirc] \\
&=[\Cbullet, (p,p),(q+2p,q+2p), \Cbullet].  
\end{align*}
}
Therefore, by the same argument as the proof of Lemma 
\ref{lem:first-Euclid}, we obtain: 
\begin{lem}\label{lem:second-Euclid}
%
For any primitive element $a$ of $H_1(\Sigma_g; {\Bbb Z})$, 
by applying 
\linebreak
$\boxplus [(\epsilon_1, \delta_1), 
\cdots, (\epsilon_g, \delta_g)]$ 
(where $\epsilon_i = 0,1$, and $\delta_i = 0,1$) several times, 
$a$ is deformed to 
$\boxplus [(\epsilon_1, \delta_1), \cdots, (\epsilon_g, \delta_g)]$ 
(where $\epsilon_i = 0,1$, and $\delta_i = 0,1$) or 
$[\Ccirc, (-1,0), \Ccirc]$. 
\qed
\end{lem}
Since $T_{-a}^2 (v) = v + 2(-a,v)(-v) = 
v + 2(a,v)v = T_a^2 (v)$, we do not need to consider 
the elements $[\Ccirc, (-1,0), \Ccirc]$. 
Hence, Lemma \ref{lem:level2-congruence-generator} follows. 
\par
\begin{figure}[ht!]
\centering
\includegraphics[height=2cm]{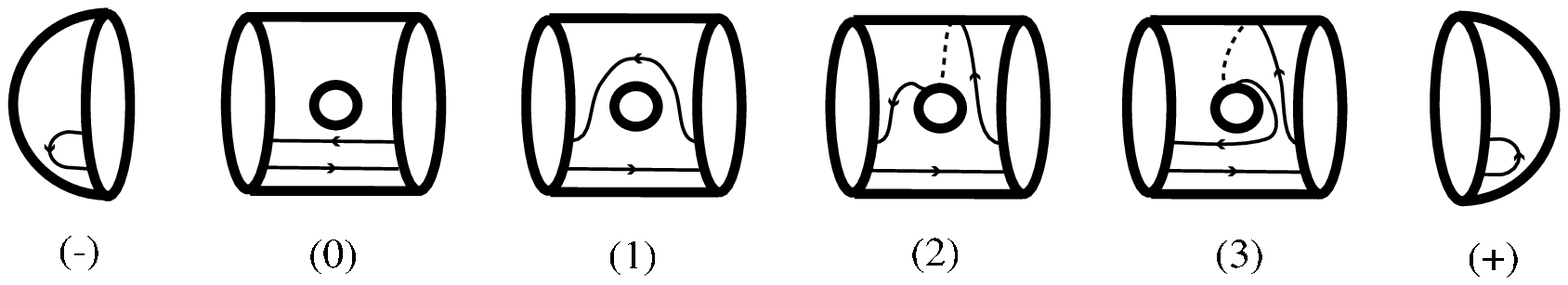}
\nocolon
\caption{}
\label{fig:realization}
\end{figure}
For each element $[(\epsilon_1, \delta_1), \cdots, (\epsilon_g, \delta_g)]$ 
(where $\epsilon_i = 0,1$, $\delta_i = 0,1$) of 
$H_1(\Sigma_g; {\Bbb Z})$, we construct an oriented simple close curve 
on $\Sigma_g$ which represent this homology class. 
For each $i$-th block, 
if $(\epsilon_i, \delta_i) = (0,0)$, we prepare (0) 
of Figure \ref{fig:realization}, 
if $(\epsilon_i, \delta_i) = (0,1)$, we prepare (1) 
of Figure \ref{fig:realization}, 
if $(\epsilon_i, \delta_i) = (1,1)$, we prepare (2) 
of Figure \ref{fig:realization}, 
if $(\epsilon_i, \delta_i) = (1,0)$, we prepare (3) 
of Figure \ref{fig:realization}. 
After that, we glue them along the boundaries and 
cap the left boundary component by (-) of Figure \ref{fig:realization} 
and the right boundary component by (+) of Figure \ref{fig:realization}. 
We denote this oriented simple closed curve on $\Sigma_g$ by 
$\{ (\epsilon_1, \delta_1), \cdots, (\epsilon_g, \delta_g) \}$. 
Here, we remark that 
the action of 
$T_{\{ (\epsilon_1, \delta_1), \cdots, (\epsilon_g, \delta_g)
\}}$ on $H_1(\Sigma_g; {\Bbb Z})$ equals  
$T_{[(\epsilon_1, \delta_1), \cdots, (\epsilon_g, \delta_g)]}$, 
and, for any $\phi$ of ${\cal M}_g$, 
$\phi \circ T_{\{ (\epsilon_1, \delta_1), \cdots, (\epsilon_g, \delta_g)
\}} \circ \phi^{-1}$ $=$ 
$T_{\phi (\{ (\epsilon_1, \delta_1), \cdots, (\epsilon_g, \delta_g)
\})}$. 
\begin{lem}\label{lem:level2-Gg-action}
%
For any $\{ (\epsilon_1, \delta_1), \cdots, (\epsilon_g, \delta_g) \}$,     
there is an element $\phi$ of $G_g$ such that 
\begin{equation}
\begin{align*}
\phi(\{ (\epsilon_1, \delta_1), \cdots, (\epsilon_g, \delta_g)\}) 
&= \{(0,1),(0,0),(0,0), \cdots, (0,0)\} \\
\text{ or } &= \{(1,1),(0,0),(0,0),\cdots,(0,0)\} \\
\text{ or } &= \{(0,0),(1,1),(0,0),\cdots,(0,0)\}. 
\end{align*}
\end{equation}
\end{lem}
\begin{figure}[ht!]
\centering
\includegraphics[height=10.5cm]{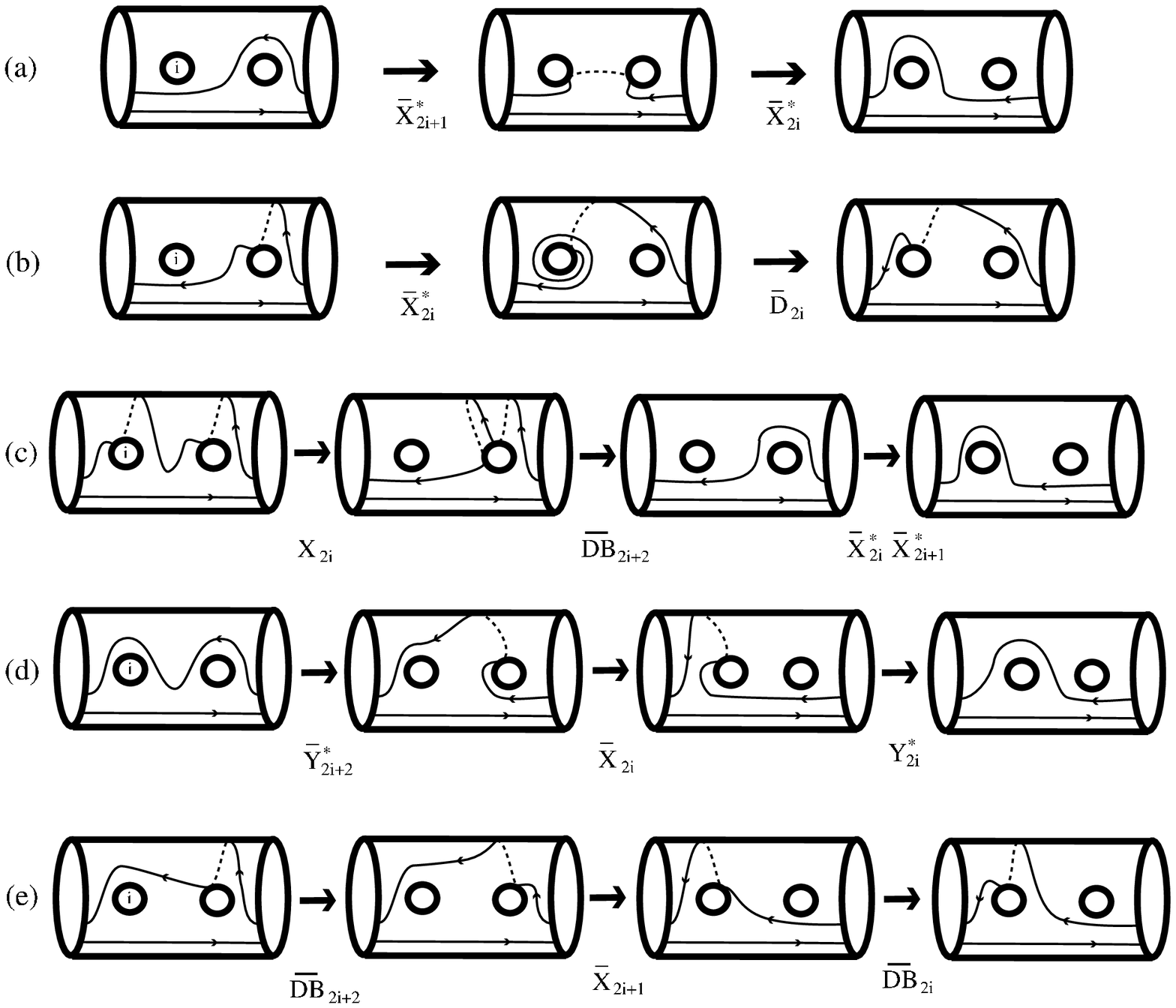}
\nocolon
\caption{}
\label{fig:elementary}
\end{figure}
\begin{proof}
If the $i$-th block is (3), by the action of $\inv{Y_{2i}}$, 
this block is changed to (1). 
Therefore, it suffices to show this lemma in the case when each block 
is not (3). 
First we investigate actions of elements of $G_g$ on adjacent 
blocks, say the $i$-th block and the $i+1$-st block. 
Each picture of Figure \ref{fig:elementary} shows the action of $G_g$ on this 
adjacent blocks. 
\begin{equation}
\begin{align*}
&(a) \text{ shows } \{ \Cbullet, (0,0),(0,1), \Cbullet \} \Ggeq 
\{ \Cbullet, (0,1), (0,0), \Cbullet \}, \\ 
&(b) \text{ shows } \{ \Cbullet, (0,0),(1,1), \Cbullet \} \Ggeq 
\{ \Cbullet, (1,1), (0,1), \Cbullet \}, \\
&(c) \text{ shows } \{ \Cbullet, (1,1),(1,1), \Cbullet \} \Ggeq 
\{ \Cbullet, (0,1), (0,0), \Cbullet \}, \\
&(d) \text{ shows } \{ \Cbullet, (0,1),(0,1), \Cbullet \} \Ggeq 
\{ \Cbullet, (0,1), (0,0), \Cbullet \}, \\
&(e) \text{ shows } \{ \Cbullet, (0,1),(1,1), \Cbullet \} \Ggeq 
\{ \Cbullet, (1,1), (0,0), \Cbullet \}.
\end{align*}
\end{equation}
For an oriented simple closed curve $x=$ 
$\{ (\epsilon_1, \delta_1), \cdots, (\epsilon_g, \delta_g)\}$, 
each of whose block is $(0,0)$ or $(0,1)$ or $(1,1)$, 
let the right most non-$(0,0)$ block be the $j$-th block. 
By the induction on $j$, we show that $x$ is $G_g$-equivalent to 
$\{(0,1),(0,0), 
$ $
(0,0), \cdots, (0,0)\}$ or 
$\{(1,1),(0,0),(0,0), \cdots, (0,0)\}$ or 
$\{(0,0),(1,1),(0,0), \cdots, 
$ $
(0,0)\}$. 
If $j$ $=$ $1$, it is trivial. 
\newline
{\it When the $j$-th block is $(0,1)$.\/} 
If each block between the first block and the $(j-1)$-st block 
is $(0,0)$, then, by repeated application of (a), 
$x$ is $G_g$-equivalent to 
$\{ (0,1), (0,0), \cdots, (0,0) \}$. 
If there is a block between the first block and the $(j-1)$-st block which
is  not $(0,0)$, by the induction hypothesis, 
the sequence from the first block to the $(j-1)$-st block 
is $G_g$-equivalent to $(0,1),(0,0),(0,0), \cdots, 
$ $
(0,0)$ or 
$(1,1),(0,0),(0,0), \cdots, (0,0)$ or $(0,0),(1,1),(0,0), \cdots, (0,0)$. 
In the first case, 
\begin{equation}
\begin{align*}
x &\Ggeq \{(0,1), (0,0), (0,0), \cdots, (0,0), (0,1),\cdots, (0,0) \}
(\text{ by the hypothesis }) \\
&\Ggeq \{(0,1), (0,1), (0,0), \cdots, (0,0),(0,0),\cdots, (0,0) \}
(\text{ by (a) }) \\
&\Ggeq \{(0,1), (0,0), (0,0), \cdots, (0,0),(0,0),\cdots, (0,0) \}
(\text{ by (d) }). 
\end{align*}
\end{equation}
In the second case, 
\begin{equation}
\begin{align*}
x &\Ggeq \{(1,1), (0,0), (0,0), \cdots, (0,0), (0,1), \cdots, (0,0) \}
(\text{ by the hypothesis }) \\
&\Ggeq \{(1,1), (0,1), (0,0), \cdots, (0,0), (0,0), \cdots, (0,0) \}
(\text{ by (a) }) \\
&\Ggeq \{(0,0), (1,1), (0,0), \cdots, (0,0), (0,0), \cdots, (0,0) \}
(\text{ by (b) }).
\end{align*}
\end{equation}
In the third case, 
\begin{equation}
\begin{align*}
x &\Ggeq \{(0,0), (1,1), (0,0), \cdots, (0,0), (0,1), \cdots, (0,0) \}
(\text{ by the hypothesis }) \\
&\Ggeq \{(0,0), (1,1), (0,1), \cdots, (0,0), (0,0), \cdots, (0,0) \}
(\text{ by (a) }) \\
&\Ggeq \{(1,1), (0,1), (0,1), \cdots, (0,0), (0,0), \cdots, (0,0) \}
(\text{ by (b) }) \\
&\Ggeq \{(1,1), (0,1), (0,0), \cdots, (0,0), (0,0), \cdots, (0,0) \}
(\text{ by (d) }) \\
&\Ggeq \{(0,0), (1,1), (0,0), \cdots, (0,0), (0,0), \cdots, (0,0) \}
(\text{ by (b) }). 
\end{align*}
\end{equation}
\newline
{\it When the $j$-th block is $(1,1)$.\/} 
If every block between the first block and $(j-1)$-st block is $(0,0)$,
then, 
\begin{equation}
\begin{align*}
x &\Ggeq \{(1,1),(0,1),(0,1)\cdots, (0,1),\cdots, (0,0)\} 
(\text{ by (b) })\\
&\Ggeq \{(1,1),(0,1),(0,0),\cdots, (0,0),\cdots, (0,0)\}  
(\text{ by (d) })\\
&\Ggeq \{(0,0),(1,1),(0,0),\cdots, (0,0),\cdots, (0,0)\}  
(\text{ by (b) }) 
\end{align*}
\end{equation}
If there is a block between the first block and the $(j-1)$-st block which
is  not $(0,0)$, by the induction hypothesis, 
the sequence from the first block to the $(j-1)$-st block 
is $G_g$-equivalent to $(0,1),(0,0),(0,0), \cdots, (0,0)$ or 
$(1,1),(0,0),(0,0), \cdots, (0,0)$ or 
$(0,0),(1,1),(0,0), \cdots, (0,0)$. 
In the first case, 
\begin{equation}
\begin{align*}
x &\Ggeq \{(0,1),(0,0),(0,0),(0,0), \cdots, (0,0),(1,1),\cdots, (0,0)\} 
(\text{ by the hypothesis })\\
&\Ggeq \{(0,1),(1,1),(0,1),(0,1), \cdots, (0,1),(0,1), \cdots, (0,0)\} 
(\text{ by (b) })\\
&\Ggeq \{(0,1),(1,1),(0,1), (0,0), \cdots, (0,0),(0,0), \cdots, (0,0)\} 
(\text{ by (d) })\\
&\Ggeq \{(1,1),(0,0),(0,1), (0,0), \cdots, (0,0),(0,0), \cdots, (0,0)\} 
(\text{ by (e) })\\
&\Ggeq \{(1,1),(0,1),(0,0), (0,0), \cdots, (0,0),(0,0), \cdots, (0,0)\} 
(\text{ by (a) })\\
&\Ggeq \{(0,0),(1,1),(0,0), (0,0), \cdots, (0,0),(0,0), \cdots, (0,0)\} 
(\text{ by (b) }). 
\end{align*}
\end{equation}
In the second case, 
\begin{equation}
\begin{align*}
x &\Ggeq \{(1,1),(0,0),(0,0),(0,0), \cdots, (0,0),(1,1), \cdots, (0,0)\} 
(\text{ by the hypothesis }) \\
&\Ggeq \{(1,1),(1,1),(0,1),(0,1), \cdots, (0,1),(0,1), \cdots, (0,0)\} 
(\text{ by (b) }) \\
&\Ggeq \{(1,1),(1,1),(0,1),(0,0), \cdots, (0,0),(0,0), \cdots, (0,0)\} 
(\text{ by (d) }) \\
&\Ggeq \{(0,1),(0,0),(0,1),(0,0), \cdots, (0,0),(0,0), \cdots, (0,0)\} 
(\text{ by (c) }) \\
&\Ggeq \{(0,1),(0,1),(0,0),(0,0), \cdots, (0,0),(0,0), \cdots, (0,0)\} 
(\text{ by (a) }) \\
&\Ggeq \{(0,1),(0,0),(0,0),(0,0), \cdots, (0,0),(0,0), \cdots, (0,0)\} 
(\text{ by (d) }). 
\end{align*}
\end{equation}
In the third case, 
\begin{equation}
\begin{align*}
x &\Ggeq \{(0,0),(1,1),(0,0),(0,0),  \cdots, (0,0), (1,1),  
\cdots, (0,0)\} (\text{ by the hypothesis }) \\
&\Ggeq \{(0,0),(1,1),(1,1),(0,1), \cdots, (0,1), (0,1), 
\cdots, (0,0)\} (\text{ by (b) }) \\
&\Ggeq \{(0,0),(0,1),(0,0),(0,1), \cdots, (0,1), (0,1), 
\cdots, (0,0)\} (\text{ by (c) }) \\
&\Ggeq \{(0,0),(0,1),(0,0),(0,1), \cdots, (0,0), (0,0), 
\cdots, (0,0)\} (\text{ by (d) }) \\
&\Ggeq \{(0,1),(0,1),(0,0),(0,0), \cdots, (0,0), (0,0), 
\cdots, (0,0)\} (\text{ by (a) }) \\
&\Ggeq \{(0,1),(0,0),(0,0),(0,0), \cdots, (0,0), (0,0), 
\cdots, (0,0)\} (\text{ by (d) }).
\end{align*}
\end{equation}
\end{proof}
By the fact that 
$T^2_{\{ (0,1),(0,0), \cdots, (0,0) \}} = D_2$, 
$T^2_{\{ (1,1),(0,0), \cdots, (0,0) \}} = (X^*_1)^2$, 
\linebreak
$T^2_{\{ (0,0),(1,1), \cdots, (0,0) \}} = (Y^*_2)^2$, 
and Lemma \ref{lem:Torelli}, 
Lemma \ref{lem:level2-congruence} is proved. 
\subsection{The modulo $2$ orthogonal group}
In this subsection, we assume $g \geq 3$.  
As in the previous subsection, let $\Phi_2\co{\cal M}_g \to \Symp (2g, {\Bbb
Z}_2)$  be the natural homomorphism. 
Let $q \co H_1(\Sigma_g; {\Bbb Z}_2) \to {\Bbb Z}_2$ be 
the quadratic form associated with the intersection form $(,)_2$ of 
$H_1(\Sigma_g; {\Bbb Z}_2)$ which satisfies  
$q (x_i) = q (y_i) = 0$ for the basis $x_i, y_i$ of 
$H_1(\Sigma_g; {\Bbb Z}_2)$ indicated on Figure \ref{fig:basis}. 
We define $\Ortho(2g, {\Bbb Z}_2)$ $=$ 
$\{ \phi \in \Aut (H_1 (\Sigma_g; {\Bbb Z}_2)) | 
q ( \phi(x)) = q (x) \text{ for any } x \in 
H_1 (\Sigma_g; {\Bbb Z}_2) \}$, 
then ${\cal SP}_g$ $=$ $\Phi_2^{-1}(\Ortho(2g, {\Bbb Z}_2))$. 
Because of Lemma \ref{lem:level2-congruence}, if we show 
$\Phi_2(G_g) = \Ortho(2g, {\Bbb Z}_2)$, then 
$G_g = {\cal SP}_g$ follows.  
For any $z$ $\in$ $H_1 (\Sigma_g; {\Bbb Z}_2)$ such that 
$q(z)=1$, we define 
${\Bbb T}_z (x)= x + (z,x)_2 \; z$. 
Then ${\Bbb T}_z$ is an element of $\Ortho(2g, {\Bbb Z}_2)$, and 
we call this a {\it ${\Bbb Z}_2$-transvection about $z$\/}. 
Dieudonn\'e \cite{Dieudonne} showed the following theorem. 
\begin{thm}\label{thm:Dieudonne}{\rm\cite[Proposition 14 on p.42]{Dieudonne}}\qua
%
When $g \geq 3$, $\Ortho(2g, {\Bbb Z}_2)$ is generated by 
${\Bbb Z}_2$-transvections. 
\end{thm}    
Let $\Lambda_g$ be the set of $z$ of $H_1 (\Sigma_g; {\Bbb Z}_2)$ 
such that $q (z) = 1$. 
For any elements $z_1$ and $z_2$ of $\Lambda_g$, we define 
$z_1 \square z_2 = z_1 + (z_2, z_1)_2 \; z_2$. 
Here, we remark that ${\Bbb T}_{z_1}^2 = \mathrm{id}$, 
${\Bbb T}_{z_2} {\Bbb T}_{z_1} {\Bbb T}_{z_2}^{-1}$ $=$ 
${\Bbb T}_{z_1 \square z_2}$ and $z_1 \square z_2 \square z_2 = z_1$. 
We denote an element 
$\epsilon_1 x_1 + \delta_1 y_1 + \cdots + \epsilon_g x_g + \delta_g y_g$ 
of $H_1 (\Sigma_g; {\Bbb Z}_2)$ by 
$[(\epsilon_1, \delta_1), \cdots, (\epsilon_g, \delta_g)]$, 
and call each $(\epsilon_i, \delta_i)$ the $i$-th block. 
$\Lambda_g$ is a set finitely generated by the operation $\square$. 
In fact, we have:
\begin{lem}\label{lem:generator-lambda-g}
%
Under the operation $\square$, $\Lambda_g$ is generated by 
$x_i + y_i$ $(1 \leq i \leq g)$, $x_i + y_i + x_{i+1}$ $(1 \leq i \leq g-1)$, 
and $x_i+x_{i+1}+y_{i+1}$ $(1 \leq i \leq g-1)$. 
\end{lem}
\begin{proof}
For an element $[(\epsilon_1, \delta_1), \cdots, (\epsilon_g, \delta_g)]$ of 
$H_1 (\Sigma_g; {\Bbb Z}_2)$, let the $j$-th block be the right most block 
which is $(1,1)$. 
When $j \geq 3$, there exist 4 cases of the combination of the $(j-1)$-st
block  and the $j$-th block: 
$[\cdots, (1,1), (1,1), \cdots]$, $[\cdots, (0,0), (1,1), \cdots]$, 
$[\cdots, (0,1), (1,1), \cdots]$, $[\cdots, (1,0), (1,1), \cdots]$. 
In each case, we can reduce $j$ at least 1. In fact, 
\begin{equation}
\begin{align*}
&[\cdots, (1,1), (1,1), \cdots] \square (x_{j-1} + x_j + y_j) 
=[\cdots, (0,1),(0,0),\cdots], \\
&[\cdots, (0,0), (1,1), \cdots] \square (x_{j-1} + y_{j-1} + x_j) 
=[\cdots, (1,1),(0,1),\cdots], \\
&[\cdots, (0,1), (1,1), \cdots] \square (x_{j-1} + x_j + y_j) 
=[\cdots, (1,1),(0,0),\cdots], \\
&[\cdots, (1,0), (1,1), \cdots] \square (x_{j-1} + y_{j-1}) 
\square (x_{j-1} + x_j + y_j) 
=[\cdots, (1,1),(0,0),\cdots]. 
\end{align*}
\end{equation}
When $j=2$, 
since $q([(\epsilon_1, \delta_1), \cdots, (\epsilon_g, \delta_g)])$=1, 
there are 3 cases of combination of the first block and the second block: 
$[(0,0),(1,1),\cdots]$, $[(1,0),(1,1),\cdots]$, or $[(0,1),(1,1),\cdots]$. 
In each case, $j$ can be reduced to $1$. In fact, 
\begin{equation}
\begin{align*}
&[(0,0),(1,1),\cdots] \square (x_1+y_1+x_2) = [(1,1),(0,1),\cdots], \\
&[(1,0),(1,1),\cdots] \square (x_1+y_1) 
\square (x_1+x_2+y_2) = [(1,1),(0,0),\cdots], \\
&[(0,1),(1,1),\cdots] \square (x_1+x_2+y_2) = [(1,1),(0,0),\cdots]. 
\end{align*}
\end{equation}
When $j=1$, if every $i$-th ($i \geq 2$) block is $(0,0)$, then 
it is $x_1 + y_1$. If there exist at least one of 
the $i$-th ($i \geq 2$) blocks which are $(1,0)$ or $(0,1)$, then,  
\begin{equation}
\begin{align*}
&[\cdots, (0,0), \overset{i}{(1,0)}, \cdots] \square (x_{i-1} + x_i + y_i) 
=[\cdots, (1,0), (0,1), \cdots], \\
&[\cdots, (1,0), \overset{i}{(0,0)}, \cdots] \square (x_{i-1} + y_{i-1} +x_i) 
=[\cdots, (0,1), (1,0), \cdots], \\
&[\cdots, (0,0), \overset{i}{(0,1)}, \cdots] \square (x_{i-1} + x_i + y_i) 
=[\cdots, (1,0), (1,0), \cdots], \\
&[\cdots, (0,1), \overset{i}{(0,0)}, \cdots] \square (x_{i-1} + y_{i-1}+x_i ) 
=[\cdots, (1,0), (1,0), \cdots].
\end{align*}
\end{equation}
Therefore, we can alter this to an element, 
each $i$-th ($i \geq 2$) block of which is $(1,0)$ or $(0,1)$. 
If the $i$-th block of this is $(0,1)$, then 
$$
[\cdots, (0,1), \cdots] \square (x_i + y_i) = [\cdots, (1,0), \cdots].
$$ 
Therefore, it suffices to consider the case when 
the first block is $(1,1)$ and other blocks are $(1,0)$. 
In this case, 
$$
[\cdots, (1,0), (1,0)] \square (x_{g-1}+y_{g-1}+x_g) 
\square (x_{g-1} + y_{g-1}) = [\cdots, (1,0), (0,0)]. 
$$
By applying the same operation repeatedly, we get 
$[(1,1),(1,0),\Ccirc]$ as a result. 
\end{proof}
This lemma and Theorem \ref{thm:Dieudonne} show:  
\begin{cor}\label{cor:generator-o-2g}
%
$\Ortho(2g, {\Bbb Z}_2)$ is generated by 
${\Bbb T}_{x_i + y_i}$ $(1 \leq i \leq g)$, 
${\Bbb T}_{x_i + y_i +x_{i+1}}$ $(1 \leq i \leq g-1)$, and 
${\Bbb T}_{x_i + x_{i+1} + y_{i+1}}$ $(1 \leq i \leq g-1)$.
\qed
\end{cor}
Since $G_g$ is a subgroup of ${\cal SP}_g$,  
$\Phi_2 (G_g) \subset \Ortho(2g, {\Bbb Z}_2)$. 
On the other hand, the fact that 
$\Phi_2(X_{2i}) = {\Bbb T}_{x_i + y_i + x_{i+1}}$ ($1 \leq i \leq g-1$), 
$\Phi_2(X_{2i+1}) = {\Bbb T}_{x_i + x_{i+1} + y_{i+1}}$ 
($1 \leq i \leq g-1$), $\Phi_2 (X_1) = {\Bbb T}_{x_1+y_1}$, 
$\Phi_2(Y_{2j}) = {\Bbb T}_{x_j + y_j}$ ($2 \leq j \leq g-1$), 
$\Phi_2(X_{2g}) = {\Bbb T}_{x_g + y_g}$, 
and Corollary \ref{cor:generator-o-2g}, 
show $\Phi_2 (G_g) \supset \Ortho(2g, {\Bbb Z}_2)$. 
Therefore we proved that, if $g \geq 3$, then ${\cal SP}_g = G_g$. 
\subsection{Genus 2 case: Reidemeister-Schreier method}
Birman and Hilden showed the following Theorem. 
\begin{thm}\label{thm:genus2}{\rm\cite{Birman-Hilden}}\qua
%
${\cal M}_2$ is generated by $C_1, C_2, C_3, C_4, C_5$ and its defining 
relations are: 
\par\noindent
{\rm(1)}\qua $C_i C_j = C_j C_i$, if $|i-j| \geq 2$, $i,j=1,2,3,4,5$, 
\par\noindent
{\rm(2)}\qua $C_i C_{i+1} C_i = C_{i+1} C_i C_{i+1}$, $i=1,2,3,4$, 
\par\noindent
{\rm(3)}\qua $(C_1 C_2 C_3 C_4 C_5)^6 = 1$, 
\par\noindent
{\rm(4)}\qua $(C_1 C_2 C_3 C_4 C_5 C_5 C_4 C_3 C_2 C_1)^2 = 1$, 
\par\noindent
{\rm(5)}\qua $C_1 C_2 C_3 C_4 C_5 C_5 C_4 C_3 C_2 C_1 \rightleftarrows C_i$, 
$i=1,2,3,4,5$, 
\par\noindent
where $\rightleftarrows$ means "commute with". 
\end{thm}
We call (1) (2) of the above relations {\it braid relations\/}. 
We will use the well-known method, called 
{\it the Reidemeister--Schreier method\/} \cite[\S 2.3]{MKS}, 
to show ${\cal SP}_2 \subset G_2$. 
We review (a part of) this method. 
\par
Let $G$ be a group generated by finite elements $g_1, \ldots, g_m$ and 
$H$ be a finite index subgroup of $G$. 
For two elements $a$, $b$ of $G$, we write $a \equiv b$ mod $H$ if 
there is an element $h$ of $H$ such that $a = hb$. 
A finite subset $S$ of $G$ is called {\it a coset representative system\/} 
for $G$ mod $H$, if, for each elements $g$ of $G$, there is only one element 
$\repre{g} \in S$ such that $g \equiv \repre{g}$ mod $H$. 
The set $\{ s g_i \repre{s g_i}^{-1} \ | \ i=1,\ldots,m,\ s \in S \}$ 
generates $H$. 
\par
\begin{figure}[ht!]
\centering
\includegraphics[height=9cm]{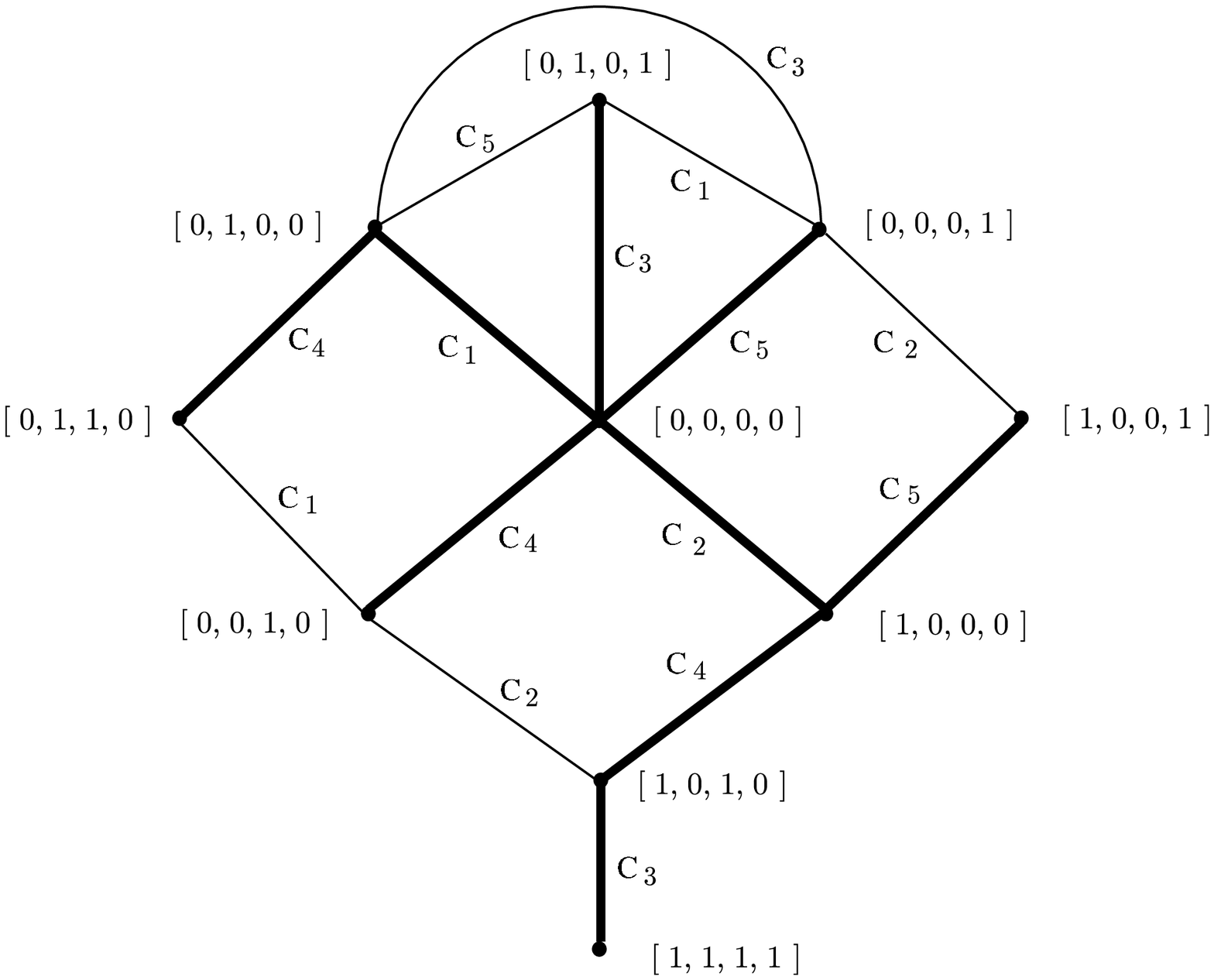}
\nocolon
\caption{}
\label{fig:action}
\end{figure}
For the sake of giving a coset representative system for ${\cal M}_2$ modulo 
${\cal SP}_2$, we will draw a graph $\Gamma$ which represents the action of 
${\cal M}_2$ on the quadratic forms of $H_1(\Sigma_2; {\Bbb Z}_2)$ with 
Arf invariants $0$. 
Let $[\epsilon_1, \epsilon_2, \epsilon_3, \epsilon_4]$ denote the quadratic form 
$q'$ of $H_1(\Sigma_2; {\Bbb Z}_2)$ such that $q'(x_1)=\epsilon_1$, 
$q'(y_1)=\epsilon_2$, $q'(x_2)=\epsilon_3$, $q'(y_2)=\epsilon_4$. 
Each vertex of $\Gamma$ corresponds to a quadratic form. 
For each generator $C_i$ of ${\cal M}_2$, we denote its action on 
$H_1(\Sigma_2; {\Bbb Z}_2)$ by $(C_i)_*$. 
For the quadratic form $q'$ indicated by the symbol 
$[\epsilon_1, \epsilon_2, \epsilon_3, \epsilon_4]$, 
let $\delta_1 = q'((C_i)_* x_1)$,  $\delta_2 = q'((C_i)_* y_1)$, 
 $\delta_3 = q'((C_i)_* x_2)$,  and $\delta_4 = q'((C_i)_* y_2)$. 
Then, we connect two vertices, corresponding to 
$[\epsilon_1, \epsilon_2, \epsilon_3, \epsilon_4]$, 
$[\delta_1, \delta_2, \delta_3, \delta_4]$ respectively, 
by the edge with the letter $C_i$. 
We remark that this action is a right action. 
For simplicity, we omit the edge whose ends are the same vertex. 
As a result, we get a graph $\Gamma$ as in Figure \ref{fig:action}. 
(Remark: The same graph was in \cite[Proof of Lemma 3.1]{Harer2}. )
In Figure \ref{fig:action}, the bold edges form a maximal tree $T$ of 
$\Gamma$. 
The words $S=\{1,\ C_1,\ C_2,\ C_3,\ C_4,\ C_5,\ C_1 C_4,\ C_2 C_4,\ 
C_2 C_5,\ C_2 C_4 C_3 \}$, which correspond to the edge paths beginning 
from $[0,0,0,0]$ on $T$, define a coset representative system for 
${\cal M}_2$ modulo ${\cal SP}_2$. 
For each element $g$ of ${\cal M}_2$, we can give a $\repre{g}$ $\in S$ with 
using this graph. 
For example, say $g=C_2 C_4 C_5 C_2$, we follow an edge path assigned 
to this word which begins from $[0,0,0,0]$, 
(note that we read words from left to right) then we arrive at the vertex 
$[0,0,1,0]$. 
The edge path on $T$ which begins from $[0,0,0,0]$ and ends at $[0,0,1,0]$ 
is $C_4$. Hence, $\repre{C_2 C_4 C_5 C_2} = C_4$. 
We list in Table \ref{tab:generators} the set of generators 
$\{ s C_i \repre{s C_i}^{-1} \ | \ i=1,\ldots,5,\ s \in S\}$ of 
${\cal SP}_g$. 
In Table \ref{tab:generators}, 
vertical direction is a coset representative system 
$S$, horizontal direction is a set of generators 
$\{C_1,\ C_2,\ C_3,\ C_4,\ C_5\}$. 
\begin{table}
\caption{Generators of ${\cal SP}_2$}
\label{tab:generators}
\begin{center}\vspace{0.1in}
\begin{tabular}{c|c c c c c}
 & $C_1$ & $C_2$ & $C_3$ & $C_4$ & $C_5$ \\
\hline
$1$ & $1$ & $1$ & $1$ & $1$ & $1$ \\
$C_1$ & $D_1$ & $X_1^*$ & $T D_5^{-1}$ & $1$ & $T D_3^{-1}$ \\
$C_2$ & $X_1$ & $D_2$ & $X_2^*$ & $1$ & $1$ \\
$C_3$ & $T D_5^{-1}$ & $X_2$ & $D_3$ & $X_3^*$ & $T D_1^{-1}$ \\ 
$C_4$ & $1$ & $1$ & $X_3$ & $D_4$ & $X_4^*$ \\
$C_5$ & $T D_3^{-1}$ & $1$ & $T D_1^{-1}$ & $X_4$ & $D_5$ \\ 
$C_1 C_4$ & $D_1$ & $X_1^*$ & $X_3$ & $D_4$ & $X_4^*$ \\
$C_2 C_4$ & $X_1$ & $D_2$ & $1$ & $D_4$ & $X_4^*$ \\
$C_2 C_5$ & $X_1$ & $D_2$ & $X_2^*$ & $X_4$ & $D_5$ \\
$C_2 C_4 C_3$ & $X_1$ & $X_3$ & $(X_2^*)^{-1} D_4 X_2^*$ & 
$X_2^*$ & $X_4^*$ \\   
\end{tabular}
\end{center}
\end{table} 
We can check this table by Figure \ref{fig:action} and braid relations. 
For example, 
\begin{align*}
&C_2 C_4 C_3 \cdot C_1 \repre{ C_2 C_4 C_3 \cdot C_1 }^{-1} = 
 C_2 C_4 C_3 C_1 (C_2 C_4 C_3)^{-1}\\
&  =C_2 C_4 C_3 C_1 C_3^{-1} C_4^{-1} C_2^{-1} = C_2 C_1 C_2^{-1} =X_1 . 
\end{align*}
This table shows that 
${\cal SP}_2 \subset G_2$ . 


\Addressesr
\end{document}

%% file: agtout.tex
\def\ifplaintex{\expandafter\ifx\csname documentclass\endcsname\relax}

\def\gtp{{\mathsurround=0pt\it $\cal G\mskip-2mu$eometry \&\ 
$\cal T\!\!$opology $\cal P\!$ublications}}  

\def\Addressesr{\bigskip
{\small \parskip 0pt \leftskip 0pt \rightskip 0pt plus 1fil \def\\{\par}
\sl\theaddress\par
\medskip
\rm Email:\stdspace\tt\theemail\hfill\rm Received:\qua\receiveddate \par}}

\def\recd{{\small Received:\qua\receiveddate\ifx\reviseddate\relax
\else\qquad Revised:\qua\reviseddate\fi\par}} 

\def\lognumber#1{\def\thelognumber{#1}}
\def\volumenumber#1{\def\thevolumenumber{#1}}
\def\volumeyear#1{\def\thevolumeyear{#1}}
\def\papernumber#1{\def\thepapernumber{#1}}
\def\pagenumbers#1#2{\def\startpage{#1}\def\finishpage{#2}}
\def\published#1{\def\publishdate{#1}}

\def\received#1{\def\receiveddate{#1}}

\def\accepted#1{\def\accepteddate{#1}}

\long\def\asciiabstract#1{\long\def\theasciiabstract{#1}}

\let\\\par\let\thelognumber\relax\let\thevolumenumber\relax
\let\thepapernumber\relax\let\thevolumeyear\relax\let\startpage\relax
\let\finishpage\relax\let\publishdate\relax\let\receiveddate\relax
\let\reviseddate\relax\let\accepteddate\relax\let\theasciititle\relax
\let\theasciiauthors\relax
\let\theasciiabstract\relax

\let\theasciiemail\relax

\ifplaintex
\font\logobig=cmssbx10 scaled 3836
\font\logomed=cmssbx10 scaled 2557
\else
\font\logobig=cmssbx10 scaled 4200
\font\logomed=cmssbx10 scaled 2800
\fi

\long\def\makeagttitle{   
\count0=\startpage
\agt\hfill      
\hbox to 45truept{\vbox to 0pt{\vglue -13truept{\logomed A\kern -.37em{\logobig 
T}\kern -.38em G}\vss}\hss}
\break
{\small Volume \thevolumenumber\ (\thevolumeyear)
\startpage--\finishpage\nl
Published: \publishdate}

\vglue .25truein

{\parskip=0pt\leftskip 0pt plus
1fil\def\\{\par\smallskip}{\Large\bf\thetitle}\par\medskip} \vglue
0.05truein

{\parskip=0pt\leftskip 0pt plus 1fil\def\\{\par}{\sc\theauthors}
\par\medskip}%
 
\vglue 0.03truein 

{\small\leftskip 25truept\rightskip 25truept{\bf Abstract}\stdspace\theabstract

{\bf AMS Classification}\stdspace\theprimaryclass
\ifx\thesecondaryclass\relax\else; \thesecondaryclass\fi\par
{\bf Keywords}\stdspace \thekeywords\par}\vglue 7truept

}   

\ifplaintex
\font\phead=cmsl9 scaled 950
\font\pnum=cmbx10 scaled 913
\font\pfoot=cmsl9 scaled 950
\headline{\vbox to 0pt{\vskip -4.5mm\line{\small\phead\ifnum
\count0=\startpage ISSN 1472-2739 (on-line) 1472-2747 (printed)
\hfill {\pnum\folio}\else\ifodd\count0\def\\{ }%
\ifx\theshorttitle\relax\thetitle\else\theshorttitle\fi\hfill{\pnum\folio}
\else\def\\{ and }{\pnum\folio}\hfill\ifx\theshortauthors\relax\theauthors
\else\theshortauthors\fi\fi\fi}\vss}}
\footline{\vbox to 0pt{\vglue 0mm\line{\small\pfoot\ifnum\count0=\startpage
\copyright\ \gtp\hfill\else
\agt, Volume \thevolumenumber\ (\thevolumeyear)\hfill\fi}\vss}}
\else
\font=cmsl9 scaled 1050
\font\lnum=cmbx10 
\font=cmsl9 scaled 1050
\makeatletter
\def\@oddhead{{\small\lhead\ifnum\count0=\startpage ISSN 1472-2739 
(on-line) 1472-2747 (printed)\hfill {\lnum\number\count0}\else\ifodd\count0
\def\\{ }\ifx\theshorttitle\relax \thetitle \else\theshorttitle\fi\hfill
{\lnum\number\count0}\else\def\\{ and }{\lnum\number\count0}
\hfill\ifx\theshortauthors\relax 
\theauthors\else\theshortauthors\fi\fi\fi}}\def\@evenhead{\@oddhead}
\def\@oddfoot{\small\lfoot\ifnum\count0=\startpage\copyright\ \gtp\hfill\else
\agt, Volume \thevolumenumber\ (\thevolumeyear)\hfill\fi}
\def\@evenfoot{\@oddfoot}
\makeatother
\fi
\let\maketitlepage\makeagttitle
\let\makeshorttitle\maketitlepage
\let\maketitle\maketitlepage

\newwrite\gtoutfile
\long\gdef\makeheadfile{  
{\def\\{, }\def\s{ }
\immediate\openout\gtoutfile head.xxx
\immediate\write\gtoutfile{To: math@arxiv.org}
\immediate\write\gtoutfile{Subject: put OR rep NNNNN:ppppp}
\immediate\write\gtoutfile{--text follows this line--}
\immediate\write\gtoutfile{Proxy-for: \ifx\theasciiauthors\relax
\theauthors\else\theasciiauthors\fi\s<\ifx\theasciiemail\relax\theemail\else\theasciiemail\fi>}
\immediate\write\gtoutfile{\noexpand\\}
\immediate\write\gtoutfile{Authors: \ifx\theasciiauthors\relax
\theauthors\else\theasciiauthors\fi}
{\def\\{ }\immediate\write\gtoutfile{Title: \ifx\theasciititle\relax
\thetitle\else\theasciititle\fi}}
\immediate\write\gtoutfile{Subj-class: GT or SG, GR etc}
\immediate\write\gtoutfile{MSC-class: \theprimaryclass\ifx\thesecondaryclass\relax\else, \thesecondaryclass\fi}
\immediate\write\gtoutfile{Journal-ref: Algebr. Geom. Topol. \thevolumenumber\s
(\thevolumeyear) \startpage-\finishpage}
\immediate\write\gtoutfile{Comments: Published by Algebraic and
Geometric Topology at}
\immediate\write\gtoutfile{\s\s\s  http://www.maths.warwick.ac.uk/agt/AGTVol\thevolumenumber/agt-\thevolumenumber-\thepapernumber.abs.html}
\immediate\write\gtoutfile{\noexpand\\}
\immediate\write\gtoutfile{}
\ifx\theasciiabstract\relax
\immediate\write\gtoutfile{\theabstract}\else
\immediate\write\gtoutfile{\theasciiabstract}\fi
\immediate\write\gtoutfile{}
\immediate\write\gtoutfile{\noexpand\\}
\immediate\write\gtoutfile{}
\immediate\closeout\gtoutfile}}  

\def\maketitlepage{\makeagttitle\makeheadfile}
\let\makeshorttitle\maketitlepage
\let\maketitle\maketitlepage

\def\ifplaintex{\expandafter\ifx\csname documentclass\endcsname\relax}

\def\gtp{{\mathsurround=0pt\it $\cal G\mskip-2mu$eometry \&\ 
$\cal T\!\!$opology $\cal P\!$ublications}}  

\def\Addressesr{\bigskip
{\small \parskip 0pt \leftskip 0pt \rightskip 0pt plus 1fil \def\\{\par}
\sl\theaddress\par
\medskip
\rm Email:\stdspace\tt\theemail\hfill\rm Received:\qua\receiveddate \par}}

\def\recd{{\small Received:\qua\receiveddate\ifx\reviseddate\relax
\else\qquad Revised:\qua\reviseddate\fi\par}} 

\def\lognumber#1{\def\thelognumber{#1}}
\def\volumenumber#1{\def\thevolumenumber{#1}}
\def\volumeyear#1{\def\thevolumeyear{#1}}
\def\papernumber#1{\def\thepapernumber{#1}}
\def\pagenumbers#1#2{\def\startpage{#1}\def\finishpage{#2}}
\def\published#1{\def\publishdate{#1}}

\def\received#1{\def\receiveddate{#1}}

\def\accepted#1{\def\accepteddate{#1}}

\long\def\asciiabstract#1{\long\def\theasciiabstract{#1}}

\let\\\par\let\thelognumber\relax\let\thevolumenumber\relax
\let\thepapernumber\relax\let\thevolumeyear\relax\let\startpage\relax
\let\finishpage\relax\let\publishdate\relax\let\receiveddate\relax
\let\reviseddate\relax\let\accepteddate\relax\let\theasciititle\relax
\let\theasciiauthors\relax
\let\theasciiabstract\relax

\let\theasciiemail\relax

\ifplaintex
\font\logobig=cmssbx10 scaled 3836
\font\logomed=cmssbx10 scaled 2557
\else
\font\logobig=cmssbx10 scaled 4200
\font\logomed=cmssbx10 scaled 2800
\fi

\long\def\makeagttitle{   
\count0=\startpage
\agt\hfill      
\hbox to 45truept{\vbox to 0pt{\vglue -13truept{\logomed A\kern -.37em{\logobig 
T}\kern -.38em G}\vss}\hss}
\break
{\small Volume \thevolumenumber\ (\thevolumeyear)
\startpage--\finishpage\nl
Published: \publishdate}

\vglue .25truein

{\parskip=0pt\leftskip 0pt plus
1fil\def\\{\par\smallskip}{\Large\bf\thetitle}\par\medskip} \vglue
0.05truein

{\parskip=0pt\leftskip 0pt plus 1fil\def\\{\par}{\sc\theauthors}
\par\medskip}%
 
\vglue 0.03truein 

{\small\leftskip 25truept\rightskip 25truept{\bf Abstract}\stdspace\theabstract

{\bf AMS Classification}\stdspace\theprimaryclass
\ifx\thesecondaryclass\relax\else; \thesecondaryclass\fi\par
{\bf Keywords}\stdspace \thekeywords\par}\vglue 7truept

}   

\ifplaintex
\font\phead=cmsl9 scaled 950
\font\pnum=cmbx10 scaled 913
\font\pfoot=cmsl9 scaled 950
\headline{\vbox to 0pt{\vskip -4.5mm\line{\small\phead\ifnum
\count0=\startpage ISSN 1472-2739 (on-line) 1472-2747 (printed)
\hfill {\pnum\folio}\else\ifodd\count0\def\\{ }%
\ifx\theshorttitle\relax\thetitle\else\theshorttitle\fi\hfill{\pnum\folio}
\else\def\\{ and }{\pnum\folio}\hfill\ifx\theshortauthors\relax\theauthors
\else\theshortauthors\fi\fi\fi}\vss}}
\footline{\vbox to 0pt{\vglue 0mm\line{\small\pfoot\ifnum\count0=\startpage
\copyright\ \gtp\hfill\else
\agt, Volume \thevolumenumber\ (\thevolumeyear)\hfill\fi}\vss}}
\else
\font=cmsl9 scaled 1050
\font\lnum=cmbx10 
\font=cmsl9 scaled 1050
\makeatletter
\def\@oddhead{{\small\lhead\ifnum\count0=\startpage ISSN 1472-2739 
(on-line) 1472-2747 (printed)\hfill {\lnum\number\count0}\else\ifodd\count0
\def\\{ }\ifx\theshorttitle\relax \thetitle \else\theshorttitle\fi\hfill
{\lnum\number\count0}\else\def\\{ and }{\lnum\number\count0}
\hfill\ifx\theshortauthors\relax 
\theauthors\else\theshortauthors\fi\fi\fi}}\def\@evenhead{\@oddhead}
\def\@oddfoot{\small\lfoot\ifnum\count0=\startpage\copyright\ \gtp\hfill\else
\agt, Volume \thevolumenumber\ (\thevolumeyear)\hfill\fi}
\def\@evenfoot{\@oddfoot}
\makeatother
\fi
\let\maketitlepage\makeagttitle
\let\makeshorttitle\maketitlepage
\let\maketitle\maketitlepage

\newwrite\gtoutfile
\long\gdef\makeheadfile{  
{\def\\{, }\def\s{ }
\immediate\openout\gtoutfile head.xxx
\immediate\write\gtoutfile{To: math@arxiv.org}
\immediate\write\gtoutfile{Subject: put OR rep NNNNN:ppppp}
\immediate\write\gtoutfile{--text follows this line--}
\immediate\write\gtoutfile{Proxy-for: \ifx\theasciiauthors\relax
\theauthors\else\theasciiauthors\fi\s<\ifx\theasciiemail\relax\theemail\else\theasciiemail\fi>}
\immediate\write\gtoutfile{\noexpand\\}
\immediate\write\gtoutfile{Authors: \ifx\theasciiauthors\relax
\theauthors\else\theasciiauthors\fi}
{\def\\{ }\immediate\write\gtoutfile{Title: \ifx\theasciititle\relax
\thetitle\else\theasciititle\fi}}
\immediate\write\gtoutfile{Subj-class: GT or SG, GR etc}
\immediate\write\gtoutfile{MSC-class: \theprimaryclass\ifx\thesecondaryclass\relax\else, \thesecondaryclass\fi}
\immediate\write\gtoutfile{Journal-ref: Algebr. Geom. Topol. \thevolumenumber\s
(\thevolumeyear) \startpage-\finishpage}
\immediate\write\gtoutfile{Comments: Published by Algebraic and
Geometric Topology at}
\immediate\write\gtoutfile{\s\s\s  http://www.maths.warwick.ac.uk/agt/AGTVol\thevolumenumber/agt-\thevolumenumber-\thepapernumber.abs.html}
\immediate\write\gtoutfile{\noexpand\\}
\immediate\write\gtoutfile{}
\ifx\theasciiabstract\relax
\immediate\write\gtoutfile{\theabstract}\else
\immediate\write\gtoutfile{\theasciiabstract}\fi
\immediate\write\gtoutfile{}
\immediate\write\gtoutfile{\noexpand\\}
\immediate\write\gtoutfile{}
\immediate\closeout\gtoutfile}}  

\def\maketitlepage{\makeagttitle\makeheadfile}
\let\makeshorttitle\maketitlepage
\let\maketitle\maketitlepage

\def\ifplaintex{\expandafter\ifx\csname documentclass\endcsname\relax}

\def\gtp{{\mathsurround=0pt\it $\cal G\mskip-2mu$eometry \&\ 
$\cal T\!\!$opology $\cal P\!$ublications}}  

\def\Addressesr{\bigskip
{\small \parskip 0pt \leftskip 0pt \rightskip 0pt plus 1fil \def\\{\par}
\sl\theaddress\par
\medskip
\rm Email:\stdspace\tt\theemail\hfill\rm Received:\qua\receiveddate \par}}

\def\recd{{\small Received:\qua\receiveddate\ifx\reviseddate\relax
\else\qquad Revised:\qua\reviseddate\fi\par}} 

\def\lognumber#1{\def\thelognumber{#1}}
\def\volumenumber#1{\def\thevolumenumber{#1}}
\def\volumeyear#1{\def\thevolumeyear{#1}}
\def\papernumber#1{\def\thepapernumber{#1}}
\def\pagenumbers#1#2{\def\startpage{#1}\def\finishpage{#2}}
\def\published#1{\def\publishdate{#1}}

\def\received#1{\def\receiveddate{#1}}

\def\accepted#1{\def\accepteddate{#1}}

\long\def\asciiabstract#1{\long\def\theasciiabstract{#1}}

\let\\\par\let\thelognumber\relax\let\thevolumenumber\relax
\let\thepapernumber\relax\let\thevolumeyear\relax\let\startpage\relax
\let\finishpage\relax\let\publishdate\relax\let\receiveddate\relax
\let\reviseddate\relax\let\accepteddate\relax\let\theasciititle\relax
\let\theasciiauthors\relax
\let\theasciiabstract\relax

\let\theasciiemail\relax

\ifplaintex
\font\logobig=cmssbx10 scaled 3836
\font\logomed=cmssbx10 scaled 2557
\else
\font\logobig=cmssbx10 scaled 4200
\font\logomed=cmssbx10 scaled 2800
\fi

\long\def\makeagttitle{   
\count0=\startpage
\agt\hfill      
\hbox to 45truept{\vbox to 0pt{\vglue -13truept{\logomed A\kern -.37em{\logobig 
T}\kern -.38em G}\vss}\hss}
\break
{\small Volume \thevolumenumber\ (\thevolumeyear)
\startpage--\finishpage\nl
Published: \publishdate}

\vglue .25truein

{\parskip=0pt\leftskip 0pt plus
1fil\def\\{\par\smallskip}{\Large\bf\thetitle}\par\medskip} \vglue
0.05truein

{\parskip=0pt\leftskip 0pt plus 1fil\def\\{\par}{\sc\theauthors}
\par\medskip}%
 
\vglue 0.03truein 

{\small\leftskip 25truept\rightskip 25truept{\bf Abstract}\stdspace\theabstract

{\bf AMS Classification}\stdspace\theprimaryclass
\ifx\thesecondaryclass\relax\else; \thesecondaryclass\fi\par
{\bf Keywords}\stdspace \thekeywords\par}\vglue 7truept

}   

\ifplaintex
\font\phead=cmsl9 scaled 950
\font\pnum=cmbx10 scaled 913
\font\pfoot=cmsl9 scaled 950
\headline{\vbox to 0pt{\vskip -4.5mm\line{\small\phead\ifnum
\count0=\startpage ISSN 1472-2739 (on-line) 1472-2747 (printed)
\hfill {\pnum\folio}\else\ifodd\count0\def\\{ }%
\ifx\theshorttitle\relax\thetitle\else\theshorttitle\fi\hfill{\pnum\folio}
\else\def\\{ and }{\pnum\folio}\hfill\ifx\theshortauthors\relax\theauthors
\else\theshortauthors\fi\fi\fi}\vss}}
\footline{\vbox to 0pt{\vglue 0mm\line{\small\pfoot\ifnum\count0=\startpage
\copyright\ \gtp\hfill\else
\agt, Volume \thevolumenumber\ (\thevolumeyear)\hfill\fi}\vss}}
\else
\font=cmsl9 scaled 1050
\font\lnum=cmbx10 
\font=cmsl9 scaled 1050
\makeatletter
\def\@oddhead{{\small\lhead\ifnum\count0=\startpage ISSN 1472-2739 
(on-line) 1472-2747 (printed)\hfill {\lnum\number\count0}\else\ifodd\count0
\def\\{ }\ifx\theshorttitle\relax \thetitle \else\theshorttitle\fi\hfill
{\lnum\number\count0}\else\def\\{ and }{\lnum\number\count0}
\hfill\ifx\theshortauthors\relax 
\theauthors\else\theshortauthors\fi\fi\fi}}\def\@evenhead{\@oddhead}
\def\@oddfoot{\small\lfoot\ifnum\count0=\startpage\copyright\ \gtp\hfill\else
\agt, Volume \thevolumenumber\ (\thevolumeyear)\hfill\fi}
\def\@evenfoot{\@oddfoot}
\makeatother
\fi
\let\maketitlepage\makeagttitle
\let\makeshorttitle\maketitlepage
\let\maketitle\maketitlepage

\newwrite\gtoutfile
\long\gdef\makeheadfile{  
{\def\\{, }\def\s{ }
\immediate\openout\gtoutfile head.xxx
\immediate\write\gtoutfile{To: math@arxiv.org}
\immediate\write\gtoutfile{Subject: put OR rep NNNNN:ppppp}
\immediate\write\gtoutfile{--text follows this line--}
\immediate\write\gtoutfile{Proxy-for: \ifx\theasciiauthors\relax
\theauthors\else\theasciiauthors\fi\s<\ifx\theasciiemail\relax\theemail\else\theasciiemail\fi>}
\immediate\write\gtoutfile{\noexpand\\}
\immediate\write\gtoutfile{Authors: \ifx\theasciiauthors\relax
\theauthors\else\theasciiauthors\fi}
{\def\\{ }\immediate\write\gtoutfile{Title: \ifx\theasciititle\relax
\thetitle\else\theasciititle\fi}}
\immediate\write\gtoutfile{Subj-class: GT or SG, GR etc}
\immediate\write\gtoutfile{MSC-class: \theprimaryclass\ifx\thesecondaryclass\relax\else, \thesecondaryclass\fi}
\immediate\write\gtoutfile{Journal-ref: Algebr. Geom. Topol. \thevolumenumber\s
(\thevolumeyear) \startpage-\finishpage}
\immediate\write\gtoutfile{Comments: Published by Algebraic and
Geometric Topology at}
\immediate\write\gtoutfile{\s\s\s  http://www.maths.warwick.ac.uk/agt/AGTVol\thevolumenumber/agt-\thevolumenumber-\thepapernumber.abs.html}
\immediate\write\gtoutfile{\noexpand\\}
\immediate\write\gtoutfile{}
\ifx\theasciiabstract\relax
\immediate\write\gtoutfile{\theabstract}\else
\immediate\write\gtoutfile{\theasciiabstract}\fi
\immediate\write\gtoutfile{}
\immediate\write\gtoutfile{\noexpand\\}
\immediate\write\gtoutfile{}
\immediate\closeout\gtoutfile}}  

\def\maketitlepage{\makeagttitle\makeheadfile}
\let\makeshorttitle\maketitlepage
\let\maketitle\maketitlepage

\def\ifplaintex{\expandafter\ifx\csname documentclass\endcsname\relax}

\def\gtp{{\mathsurround=0pt\it $\cal G\mskip-2mu$eometry \&\ 
$\cal T\!\!$opology $\cal P\!$ublications}}  

\def\Addressesr{\bigskip
{\small \parskip 0pt \leftskip 0pt \rightskip 0pt plus 1fil \def\\{\par}
\sl\theaddress\par
\medskip
\rm Email:\stdspace\tt\theemail\hfill\rm Received:\qua\receiveddate \par}}

\def\recd{{\small Received:\qua\receiveddate\ifx\reviseddate\relax
\else\qquad Revised:\qua\reviseddate\fi\par}} 

\def\lognumber#1{\def\thelognumber{#1}}
\def\volumenumber#1{\def\thevolumenumber{#1}}
\def\volumeyear#1{\def\thevolumeyear{#1}}
\def\papernumber#1{\def\thepapernumber{#1}}
\def\pagenumbers#1#2{\def\startpage{#1}\def\finishpage{#2}}
\def\published#1{\def\publishdate{#1}}

\def\received#1{\def\receiveddate{#1}}

\def\accepted#1{\def\accepteddate{#1}}

\long\def\asciiabstract#1{\long\def\theasciiabstract{#1}}

\let\\\par\let\thelognumber\relax\let\thevolumenumber\relax
\let\thepapernumber\relax\let\thevolumeyear\relax\let\startpage\relax
\let\finishpage\relax\let\publishdate\relax\let\receiveddate\relax
\let\reviseddate\relax\let\accepteddate\relax\let\theasciititle\relax
\let\theasciiauthors\relax
\let\theasciiabstract\relax

\let\theasciiemail\relax

\ifplaintex
\font\logobig=cmssbx10 scaled 3836
\font\logomed=cmssbx10 scaled 2557
\else
\font\logobig=cmssbx10 scaled 4200
\font\logomed=cmssbx10 scaled 2800
\fi

\long\def\makeagttitle{   
\count0=\startpage
\agt\hfill      
\hbox to 45truept{\vbox to 0pt{\vglue -13truept{\logomed A\kern -.37em{\logobig 
T}\kern -.38em G}\vss}\hss}
\break
{\small Volume \thevolumenumber\ (\thevolumeyear)
\startpage--\finishpage\nl
Published: \publishdate}

\vglue .25truein

{\parskip=0pt\leftskip 0pt plus
1fil\def\\{\par\smallskip}{\Large\bf\thetitle}\par\medskip} \vglue
0.05truein

{\parskip=0pt\leftskip 0pt plus 1fil\def\\{\par}{\sc\theauthors}
\par\medskip}%
 
\vglue 0.03truein 

{\small\leftskip 25truept\rightskip 25truept{\bf Abstract}\stdspace\theabstract

{\bf AMS Classification}\stdspace\theprimaryclass
\ifx\thesecondaryclass\relax\else; \thesecondaryclass\fi\par
{\bf Keywords}\stdspace \thekeywords\par}\vglue 7truept

}   

\ifplaintex
\font\phead=cmsl9 scaled 950
\font\pnum=cmbx10 scaled 913
\font\pfoot=cmsl9 scaled 950
\headline{\vbox to 0pt{\vskip -4.5mm\line{\small\phead\ifnum
\count0=\startpage ISSN 1472-2739 (on-line) 1472-2747 (printed)
\hfill {\pnum\folio}\else\ifodd\count0\def\\{ }%
\ifx\theshorttitle\relax\thetitle\else\theshorttitle\fi\hfill{\pnum\folio}
\else\def\\{ and }{\pnum\folio}\hfill\ifx\theshortauthors\relax\theauthors
\else\theshortauthors\fi\fi\fi}\vss}}
\footline{\vbox to 0pt{\vglue 0mm\line{\small\pfoot\ifnum\count0=\startpage
\copyright\ \gtp\hfill\else
\agt, Volume \thevolumenumber\ (\thevolumeyear)\hfill\fi}\vss}}
\else
\font=cmsl9 scaled 1050
\font\lnum=cmbx10 
\font=cmsl9 scaled 1050
\makeatletter
\def\@oddhead{{\small\lhead\ifnum\count0=\startpage ISSN 1472-2739 
(on-line) 1472-2747 (printed)\hfill {\lnum\number\count0}\else\ifodd\count0
\def\\{ }\ifx\theshorttitle\relax \thetitle \else\theshorttitle\fi\hfill
{\lnum\number\count0}\else\def\\{ and }{\lnum\number\count0}
\hfill\ifx\theshortauthors\relax 
\theauthors\else\theshortauthors\fi\fi\fi}}\def\@evenhead{\@oddhead}
\def\@oddfoot{\small\lfoot\ifnum\count0=\startpage\copyright\ \gtp\hfill\else
\agt, Volume \thevolumenumber\ (\thevolumeyear)\hfill\fi}
\def\@evenfoot{\@oddfoot}
\makeatother
\fi
\let\maketitlepage\makeagttitle
\let\makeshorttitle\maketitlepage
\let\maketitle\maketitlepage

\newwrite\gtoutfile
\long\gdef\makeheadfile{  
{\def\\{, }\def\s{ }
\immediate\openout\gtoutfile head.xxx
\immediate\write\gtoutfile{To: math@arxiv.org}
\immediate\write\gtoutfile{Subject: put OR rep NNNNN:ppppp}
\immediate\write\gtoutfile{--text follows this line--}
\immediate\write\gtoutfile{Proxy-for: \ifx\theasciiauthors\relax
\theauthors\else\theasciiauthors\fi\s<\ifx\theasciiemail\relax\theemail\else\theasciiemail\fi>}
\immediate\write\gtoutfile{\noexpand\\}
\immediate\write\gtoutfile{Authors: \ifx\theasciiauthors\relax
\theauthors\else\theasciiauthors\fi}
{\def\\{ }\immediate\write\gtoutfile{Title: \ifx\theasciititle\relax
\thetitle\else\theasciititle\fi}}
\immediate\write\gtoutfile{Subj-class: GT or SG, GR etc}
\immediate\write\gtoutfile{MSC-class: \theprimaryclass\ifx\thesecondaryclass\relax\else, \thesecondaryclass\fi}
\immediate\write\gtoutfile{Journal-ref: Algebr. Geom. Topol. \thevolumenumber\s
(\thevolumeyear) \startpage-\finishpage}
\immediate\write\gtoutfile{Comments: Published by Algebraic and
Geometric Topology at}
\immediate\write\gtoutfile{\s\s\s  http://www.maths.warwick.ac.uk/agt/AGTVol\thevolumenumber/agt-\thevolumenumber-\thepapernumber.abs.html}
\immediate\write\gtoutfile{\noexpand\\}
\immediate\write\gtoutfile{}
\ifx\theasciiabstract\relax
\immediate\write\gtoutfile{\theabstract}\else
\immediate\write\gtoutfile{\theasciiabstract}\fi
\immediate\write\gtoutfile{}
\immediate\write\gtoutfile{\noexpand\\}
\immediate\write\gtoutfile{}
\immediate\closeout\gtoutfile}}  

\def\maketitlepage{\makeagttitle\makeheadfile}
\let\makeshorttitle\maketitlepage
\let\maketitle\maketitlepage